\documentclass[10pt]{amsart2000}
\newenvironment{proofof}[1]{\medskip\noindent
               \textit{#1.}}{\hfill{$\square$}\\}

\newenvironment{EG}[1]{{\vspace{1 ex}}\noindent {\sc Example.}{#1}
                {\vspace{1 ex}}}

\input xy
 \xyoption{all}

\usepackage{latexsym}
\usepackage{amsmath2000}

\usepackage{amssymb}
\usepackage{amscd2000}
\usepackage{epsfig}
\usepackage{psfig}

\usepackage{multicol}
\usepackage{amsfonts}
\usepackage{palatino}

\newtheorem{thm}{Theorem}[section]
\newtheorem{theorem}[thm]{Theorem}

\newtheorem{lemma}[thm]{Lemma}
\newtheorem{prop}[thm]{Proposition}

\newtheorem{DEF}[thm]{Definition}
\newtheorem{remark}[thm]{Remark}

\numberwithin{equation}{section}

\def\to{\rightarrow}

\newcommand{\Z}{\mathbb Z}
\newcommand{\F}{\mathbb F}
\newcommand{\R}{\mathbb R}
\newcommand{\C}{\mathbb C}
\newcommand{\D}{{\mathcal D}_n}

\def\iso{\cong}

\def\nab{\nabla}
\def\star{\mbox{star}}

\newcommand{\st}[1]{\ensuremath{^{\scriptstyle \textrm{#1}}}}
\def\Hol{\mbox{Hol}}

\newcommand{\p}{{\mathfrak p}}
\newcommand{\q}{{\mathfrak q}}

\DeclareMathOperator{\sign}{sign}
\DeclareMathOperator{\spann}{span}

\title[How is a graph like a manifold?]{
How is a graph like a manifold?
}
\date{ \today }

\author[Bolker]{Ethan D. Bolker}
\address[Ethan D. Bolker]{Department of Mathematics and Computer
Science\\
University of Massachusetts-Boston\\
Boston, MA 02125-3393}
\email{eb@cs.umb.edu}
\author[Guillemin]{Victor W. Guillemin}
\address[Victor W. Guillemin]{Department of Mathematics\\
Massachusetts Institute of Technology\\
Cambridge, MA 02139}
\email{vwg@math.mit.edu}
\author[Holm]{Tara S. Holm}
\address[Tara Holm]{Department of Mathematics\\Massachusetts Institute
of Technology 2-488\\Cambridge MA 02139}
\email{tsh@math.mit.edu}

\begin{document}

\bibliographystyle{plain}

\begin{abstract}
In this article, we discuss some classical problems in combinatorics
which can be solved by exploiting analogues between graph theory and
the theory of manifolds.  One well-known example is the McMullen
conjecture, which was settled twenty years ago by Richard Stanley by
interpreting certain combinatorial invariants of convex polytopes as
the Betti numbers of a complex projective variety.  Another example is
the classical parallel redrawing problem, which turns out to be
closely related to the problem of 
computing the second Betti number of a complex compact
$(\C^*)^n$-manifold.
\end{abstract}

\maketitle




\section{Introduction}\label{se:intro}

Some recent developments in the theory of group actions on complex
manifolds have revealed unexpected connections between the geometry of
manifolds and the geometry of graphs.  Our purpose in this
semi-expository paper is twofold: first, to explore these connections;
and second, to discuss some problems in graph theory which ``manifold''
ideas have helped to clarify.  One such problem is that of
counting the number of $n-k$ dimensional faces of a simple
$n$-dimensional convex polytope.  This number can be expressed as a
sum
\begin{equation}\label{eq:1}
f_{n-k} = \sum_{\ell=0}^k\binom{n-\ell}{k-\ell}\beta_{\ell},
\end{equation}
where the $\beta_{\ell}$'s are positive integers.  A celebrated
conjecture of McMullen \cite{McMullen} asserts that these integers satisfy the
identities 
\begin{equation}\label{eq:2}
\beta_k=\beta_{n-k}
\end{equation}
and the inequalities
\begin{equation}\label{eq:3}
\beta_0\leq\beta_1\leq\cdots\leq\beta_r,
\end{equation}
where $r=\left\lceil\frac{n}{2}\right\rceil$. In 1980, Richard Stanley
\cite{Stanley} solved McMullen's conjecture; his proof involved
showing that the $\beta_k$'s are the Betti numbers of a complex
projective variety, and hence that (\ref{eq:2}) is just Poincar\'e
duality and (\ref{eq:3}) is the ``hard Lefshetz'' theorem.

Another such problem is the classical {\em parallel redrawing}
problem: Given a graph embedded in $\R^n$, how many ways can one
reposition the vertices so that the edges of the deformed graph are
parallel to the edges of the original graph?  We show that the number
of such deformations can be counted by a combinatorial invariant
involving the zeroth and first ``Betti numbers'' of the graph.

In the remainder of this section, we describe how graph theoretical
structures occur in the study of group actions on complex manifolds.
This section provides much of the geometric motivation for the graph
theory that follows.  Although the remainder of the paper does not
depend on this section, the reader is strongly encouraged to acquaint
herself with the ideas discussed here, particularly with the examples
given in Section~\ref{subse:EGs}. The rest of the paper is graph
theoretical.  We give purely graph-theoretic definitions of various
combinatorial notions associated with the graphs coming from
manifolds. However, we make our definitions sufficiently broad so as
to apply to many graphs which are {\em not} associated with group
actions on manifolds. In Sections~\ref{se:conn} through
\ref{se:betti}, we define the notions of {\em connection} and {\em
Morse function} on a regular graph, and using ``Morse theory,'' we
define the {\em Betti numbers} of a graph.  Then in
Section~\ref{se:simple}, we take up the problem mentioned above:
counting the number of $n-k$ dimensional faces of a simple
$n$-dimensional convex polytope, $\Delta$, and show that the
$\beta_\ell$'s in Equation~(\ref{eq:1}) are just the Betti
numbers of $\Delta$.  We also provide in this section a brief description
of Stanley's proof of the McMullen conjectures.

In Section~\ref{se:poly}, we discuss another problem 
in which graph theoretical Betti numbers play an important role.  This
problem is a polynomial interpolation problem for graphs which
generalizes, in some sense, the classical problem of Lagrangian
interpolation for polynomials of one variable.  For instance, for the
complete graph on $n$ vertices, this problem consists of finding a
polynomial in one variable,
$$
p(z)=\sum_{i=1}^{n-1}g_iz^i
$$
whose coefficients $g_i$ are polynomials in $x_1,\dots,x_n$ such that
the expressions
\begin{equation}\label{eq:4}
f_i=p(z)|_{z=x_i}\ \ i=1,\dots,n
\end{equation}
are pre-assigned polynomials in the $x_i$'s.  
We will show that \eqref{eq:4} is solvable if and only if
\begin{equation}\label{eq:5}
f_i-f_j=0\mod x_i-x_j;
\end{equation}
so \eqref{eq:4} can be regarded as a formula for the general solution
of the system of equations \eqref{eq:5}.  One encounters systems of
equations of the type \eqref{eq:5} involving more complicated graphs
in the area of combinatorics known as spline theory.  Furthermore,
they come up, somewhat unexpectedly, in the theory of manifolds.  The
computation of the equivariant cohomology rings of manifolds we
describe in the examples below involves solving such equations.  We will
not attempt to discuss equivariant cohomology in this article, but
interested readers can find details in \cite{GKM}, \cite{TW}, and
\cite{GZ2}. 

In Section~\ref{se:poly}, we will show how a connection on a graph
enables one to construct many explicit solutions to interpolation
equations of this type. In Section~\ref{se:para}, we will show that
the dimension of the space of solutions of degree $k$ is given by a
formula similar to the McMullen formula \eqref{eq:1}.  Moreover for
$k=1$, this formula will 
solve
the parallel redrawing
problem for a large class of interesting graphs.  In
Section~\ref{sec:8}, we will discuss some connections between this
formula and McMullen's formula, and show that this resemblance between
the two is not entirely fortuitous.  More explicitly, we show that the
face counting problem for polytopes can be viewed as an interpolation
problem for polynomials in ``anti-commuting'' variables.  Finally, in
Section~\ref{se:EGs}, we provide several families of examples of
graphs to which this theory applies
and some open questions.
In the interest of brevity, many
details of this section are left to the reader.

We are grateful to many colleagues and friends for enlightening
comments about various aspects of this paper.  In particular, the
proof of Theorem~\ref{th:interp} is,
with minor modifications, identical to that of Theorem~2.4.1 in
\cite{GZ2}.  We are grateful to Catalin Zara for 
letting us reproduce this argument here, and also for a number of
helpful comments and insights.  We would also like to thank Werner
Ballmann, Sara Billey, Daniel Biss, Rebecca Goldin, Allen Knutson, Sue
Tolman, David Vogan, and Walter Whitely for helpful conversations.

\subsection{GKM Manifolds}\label{subse:GKM}

Henceforth, we denote by \(\C^*\) the multiplicative group
of non-zero complex numbers, and by \(G\) the group 
\[
G=\C^*\times\cdots\times\C^*=(\C^*)^n.
\]
Let \(M\) be a \(d\)-dimensional compact complex manifold, and
\(\rho:G\times M\to M\) a holomorphic action of \(G\) on \(M\).  The
action \(\rho\) is called a {\em GKM action} if there are a finite
number of \(G\) fixed points and a finite number of \(G\) orbits of
(complex) dimension one.  

An important example is the Riemann sphere $\C P^1$ with the standard
action of \(\C^*\).  This is the action
\[
c\cdot[z_0:z_1]=[c\cdot z_0:z_1],
\]
where \([z_o:z_1]\) are standard homogeneous coordinates on \(\C
P^1\).  This action has two fixed points, the origin, \([0:1]\), and
the point at infinity, \([1:0]\).  The complement of these points is a
single \(\C^*\) orbit.

When \(\rho\) is a GKM action, we let $V=M^G=\{p_1,\dots,p_{\ell}\}$
be the set of fixed points of $\rho$ and $E=\{ e_1,\dots,e_N\}$ the
set of one-dimensional orbits.  It is not hard to show that thes
objects satisfy
\begin{enumerate}
\item\label{it1} The closure of each \(e\in E\) is an embedded copy of
the complex projective line, \(\C P^1\);

\item\label{it2} If the closure of \(e_i\) intersects the closure of
\(e_j\), the intersection is a single point \(p\).  Moreover, \(p\in
M^G\);

\item\label{it3} The closure of each \(e\) contains exactly two fixed
points;

\item\label{it4} Every fixed point is contained in the closure of
exactly \(d\) one-dimensional orbits; and

\item\label{it5} The action of \(G\) on the closure of \(e\) is
isomorphic to the standard action of \(\C^*\) on \(\C P^1\).
\end{enumerate}
This last item needs some explaining, which we postpone until the next
section.

\subsection{GKM graphs}\label{subse:GKMgr}
Goresky, Kottwitz, and MacPherson point out in \cite{GKM} that the
properties~(\ref{it1}) through (\ref{it4}) are nicely expressed in the
language of graph theory. Let $\Gamma$ be the graph
having vertex set $V=M^G$ and edge set $E$, where $e\in E$ joins $p$
and $q$ in $V$ exactly when $p$ and $q$ are the two fixed points of
$G$ in the closure of $e$.  We call $\Gamma$ the {\em GKM graph} of
the pair $(M,\rho )$.  

Properties~(\ref{it1}) through (\ref{it4}) tell us that $\Gamma$ is a
particularly nice graph.  For instance, properties~(\ref{it2}) and
(\ref{it3}) imply that there are no loops and no multiple edges in
$\Gamma$.
Property~(\ref{it4}) tells us that $\Gamma$ is a {\em $d$-regular
graph}: each vertex has valence $d$. 

Property~(\ref{it5}) does not have a purely combinatorial
interpretation. It specifies a way to assign a vector in $\R^n$ to
each edge of $\Gamma$, and so provides a kind of embedding of
$\Gamma$.  Suppose $e$ is 
the edge joining $p$ and $q$.  To specify an
orientation of $e$, we decide which of $p$ and $q$ is the initial
vertex $\iota(e)$ of $e$ and which the terminal vertex $\tau(e)$ of
$e$.  Property~(\ref{it5}) says that the action of \(G\) on
the closure of \(e\) is isomorphic to the standard action of $\C^*$ on
$\C P^1$.  This means that there is a group homomorphism
\[
\chi_e:G\to\C^*
\]
and a bi-holomorphic map
\[
f_e:\overline{e}\to\C P^1
\]
such that \(f_e(g\cdot p)=\chi_e(g)\cdot f_e(p)\) for all \(g\in G\)
and for all \(p\in\overline{e}\).  The homomorphism \(\chi_e\) and the
map \(f_e\) are more or less unique.  The only way that \(f_e\) can be
altered is by composing it with a bi-holomorphic map of \(\C P^1\)
which commutes with the action of \(\C^*\), and the group of such
bi-holomorphisms consists of \(\C^*\) itself and the involution
\[
\sigma([z_0:z_1]) = [z_1:z_0]
\]
which flips the fixed points at zero
and infinity.  Moreover, if \(f_e\) is composed with an element of
\(\C^*\), \(\chi_e\) is unchanged, and composition with \(\sigma\)
interchanges \(\chi_e\) with \(\chi_e^{-1}\).  In particular,
\(\chi_e\) is determined when we specify which of the two points of
\(\overline{e}-e\) maps to zero and which 
to infinity.

In the coordinate system determined by the expression of $G$ as a
product of $\C^*$, this intertwining homomorphism has the form
$$
\chi_e(c_1,\dots,c_n)=c_1^{\alpha_1(e)}\cdots c_n^{\alpha_n(e)}.
$$
Note that the $\alpha_r(e)$ are integers
because they are characters of the group $S^1\subset\C^*$. So to
every oriented edge $e$ of $\Gamma$, we can attach the $n$-tuple of
integers 
$$
\alpha(e)=(\alpha_1(e),\dots,\alpha_n(e))\in\Z^n.
$$
If we 
reverse the orientation of $e$, then 
$\chi_e$ 
becomes $\chi_e^{-1}$ and 
$\alpha(e)$ 
becomes $-\alpha(e)$.  We call 
the map
$$
\alpha:E\to\Z^n\subset\R^n
$$ 
the {\em axial function} of the graph
$\Gamma$.  The rationale for the name ``axial function'' is that the
map $\alpha$ describes how each of the two-spheres $e=\C P^1=S^2$ gets
rotated about its axis by the torus subgroup
$$
G_{\R}=\underbrace{S^1\times\cdots\times S^1}_{n} 
$$ 
of $G$.

\subsection{Examples}\label{subse:EGs}

We describe below a few typical examples of GKM graphs.  The
reader should make herself familiar with these examples since they
will resurface frequently in the remainder of the paper.

\begin{EG}{
{\bf The $n$-simplex.} Let $M=\C P^{n-1}$ be complex projective $n-1$
space, represented using homogeneous coordinates in $\C^{n}$,  and
let $G$ act on $M$ by the product action
$$
g\cdot [z_1:\cdots:z_n]=[g_1\cdot z_1:\cdots:g_n\cdot z_n].
$$
The fixed points of this action are the points
$[0:\cdots:0:1:0:\cdots:0]$, and the one-dimensional orbits are the
sets $e_{i,j}=\{ [0:\cdots:z:\cdots:w:\cdots:0]\}$, points with
non-zero entries in just two components.  It follows easily that the
intersection graph is $K_n$, the complete graph on $n$ points.
Moreover, if we choose zero and infinity in $\overline{e_{i,j}}$ to be
the point with $i$th coordinate zero and $j$th coordinate zero
respectively, then the intertwining homomorphism is
$$
\chi_{e_{i,j}}(c_1,\dots,c_n) = c_i^1\cdot c_j^{-1},
$$
since in homogeneous coordinates,
$$
[c_ic_j^{-1}z:w]=[c_iz:c_jw].
$$
Thus, $\alpha(e_{i,j})=(0,\dots,1,0,\dots,-1,0,\dots,0)$, with a $1$ in the
$i$th position and a $-1$ in the $j$th position.
If we embed the graph $K_n$ 
in $\R^n$ by sending the vertex
$[0:\cdots:0:1:0:\cdots:0]$ to the point $(0,\dots,0,1,0,\dots,0)$,
this computation shows that  the edge joining two vertices does indeed have the direction
$\alpha(e)$.  Thus, we have proved the following proposition.

\begin{prop}
The graph $\Gamma$ of $M$ is the one-skeleton (vertices and edges) of
the $(n-1)$-simplex
$$
\{ (x_1,\dots,x_n)\in\R^n_+\ |\ x_1+\cdots +x_n=1\},
$$
and for every edge $e$ of this embedded graph, $\alpha(e)$ is the embedded
image of $e$.
\end{prop}
}\end{EG}

\begin{EG}{
{\bf The hypercube.} Let $M$ be
\begin{equation}\label{3.1}
M=\underbrace{\C P^1\times\cdots\times\C P^1}_{n}.
\end{equation}
Then $G$ acts on $M$ by the product action,
$$
c\cdot (p_1,\dots,p_n)=(c_1\cdot p_1,\dots,c_n\cdot p_n),
$$
and the fixed points for this action are the points
$$
(\varepsilon_1,\dots,\varepsilon_n),
$$
where $\varepsilon_i\in\{ 0,\infty\}$ is zero or the point at infinity
of the $i$th copy of $\C P^1$ in the product (\ref{3.1}).  It is easy
to see that the one-dimensional orbits of $G$ are just the sets
$$
(\varepsilon_1,\dots,\varepsilon_{i-1},p_i,\varepsilon_{i+1},
\dots,\varepsilon_n),
$$
where $p_i\neq \varepsilon_i$.

We leave the proof of the following proposition as an exercise:

\begin{prop}
The graph of $M$ is the one-skeleton of the hypercube
$$
\underbrace{[0,1]\times\cdots\times [0,1]}_{n}
$$
in $\R^n$,
and for every edge $e$ of this embedded graph, $\alpha(e)$ is the embedded
image of $e$.

\end{prop}

Note that this embedding of the graph does give us the axial function
as well.  This construction naturally generalizes to the product of
any set of GKM manifolds.
}\end{EG}

\begin{EG}{
{\bf The Johnson graph.}  Projective space, $\C P^{n-1}$, is really
just the set of lines in complex $n$-space $\C^n$.  This example
generalizes that construction.  The group $G$ acts on $\C^n$ with the
product action
\begin{equation}\label{3.2}
c\cdot (z_1,\dots,z_n)=(c_1\cdot z_1,\dots,c_n\cdot z_n).
\end{equation}
Let $M$ be the $k$-Grassmannian $\mathcal{G}r(k,n)$ : the set of
all (complex) $k$-dimensional subspaces of $\C^n$.  The action of
$G$ on $\C^n$ induces an action of $G$ on $M$, and it is
easy to see that the points of $M$ which are fixed by this action,
that is the $k$-dimensional subspaces of $\C^n$ which are mapped to
themselves by (\ref{3.2}) are the subspaces
$$
\C^n_I=\{ (z_1,\dots,z_n)\ |\ z_i=0 \mbox{ for } i\not\in I\},
$$
where $I$ is a $k$-element subset of $\{ 1,\dots,n\}$.  In other
words, the fixed points for the action of $G$ on $M$ are indexed by
the $k$-element subsets of $\{ 1,\dots,n\}$.  The one-dimensional
orbits of $G$ are not much harder to describe.  Let $I$ and $J$ be
$k$-element subsets of $\{ 1,\dots,n\}$ 
whose
intersection has order $k-1$.  Then the collection of all
$k$-dimensional subspaces $V$ of $\C^n$ satisfying
\begin{equation}\label{3.3}
\C^n_I\cap\C^n_J\subset V\subset \C^n_I+\C^n_J,\, V\neq \C^n_I,\C^n_J
\end{equation}
is a one-dimensional orbit.  Moreover, the sets (\ref{3.3}) are the
only one-dimensional orbits.  Hence 

\begin{prop}
The GKM graph $\Gamma$ of $M$ is the graph whose vertices correspond
to $k$-element subsets of $\{ 1,\dots,n\}$, and two such subsets $I$
and $J$ are adjacent if and only if their intersection is a
$(k-1)$-element set.
If we think of each $k$-element set as a vector in in $\R^n$ with $k$
1's and $n-k$ 0's in the obvious way then we have embedded $\Gamma$  
in $\R^n$. The axial function at an edge $e$ is the embedded
image of $e$.
\end{prop}

This graph is known as the {\em Johnson graph}, $J(k,n)$.  
$J(2,4)$ is the 
$1$-skeleton of the octahedron, 
shown below.
\begin{figure}[h]
\centerline{
\epsfig{figure=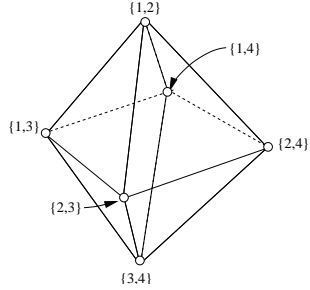,height=1.5in}
}
\centerline{
\parbox{4.5in}{\caption[$J(2,4)$]{{\small The Johnson graph
$J(2,4)$ as an octahedron.}}}
}
\end{figure}

Although $J(2,4)$ is the cross polytope
in dimension three, 
there is in general
no relation between the cross polytopes and the Johnson graphs.
}\end{EG}

\begin{EG}{
{\bf The one skeleton of a simple convex polytope}
A convex $n$-dimensional polytope $\Delta$ in $\R^n$ is
{\em simple} if 
exactly $n$ edges 
meet at each vertex.  Given such a
polytope, one can associate with it a complex $n$-dimensional space
$M_\Delta$ and an action on $M$ of $(\C^*)^n$ such that the
$k$-dimensional orbits are in one-to-one correspondence with the
$k$-dimensional faces of $\Delta$.  In particular, the zero- and
one-dimensional orbits correspond to the vertices and edges of
$\Delta$.  Thus, this action is a GKM action, and its graph is the
one-skeleton of $\Delta$.  We will describe this example in more
detail in Section~\ref{se:simple}.  The first two examples described
above, the $n$-simplex and the hypercube, are special cases of this
example, but
the polytopes associated with the Johnson
graph are not simple. }
\end{EG}

\subsection{How is a graph like a manifold?}\label{subse:graMan}

Several recent articles about the topology and
geometry of GKM manifolds exploit the fact that topological
properties of $M$ have combinatorial implications for $\Gamma$ and
that, conversely, given information about $\Gamma$, one can draw
conclusions about the topology of $M$.  See \cite{GKM}, \cite{TW},
\cite{GZ1}, \cite{GZ2}, \cite{KR}, and \cite{LLY}. Throughout these
articles, one glimpses a kind of dictionary in which manifold
concepts translate to graph concepts and vice versa.  Our modest goal
in this article is to examine some entries in the ``graph'' column of
this dictionary for their own intrinsic (graph theoretic) interest.
For instance, in Section~\ref{se:conn} we will give a combinatorial
definition of the notion of a {\em connection} on a graph $\Gamma$,
{\em geodesics}, {\em totally geodesic subgraphs} of $\Gamma$,
and the {\em holonomy group} associated to a connection.  In
Section~\ref{se:axial}, we define an {\em axial function} on a graph.
In Section~\ref{se:betti}, we discuss {\em Morse functions} and define
the combinatorial {\em Betti numbers} of $\Gamma$.  
Then we use those topological notions to prove theorems about graphs.
Finally, there are
interesting graphs for which these notions make sense
which are {\em not} GKM graphs.  
For example, the complete bipartite
graph $K_{n,n}$ has many of these combinatorial structures, but is not
associated with a manifold 
We discuss this and
other examples in more detail in Section~\ref{se:EGs}.

\section{Connections and geodesic subgraphs} \label{se:conn}

Let $\Gamma = (V,E)$ be a graph with finite vertex set $V$ and edge
set $E$. We count each edge twice, once with each of its
two possible orientations.  When $x$ and $y$ are adjacent vertices we
write $e = (x,y)$ for the edge from $x$ to $y$ and $e^{-1}=(y,x)$ for
the edge from $y$ to $x$.  Given an oriented edge $e=(x,y)$ , we write
$x=\iota(e)$ for the initial vertex and $y=\tau(e)$ for the terminal
vertex.

\begin{DEF}
The {\em star} of a vertex $x$, written $\star(x)$,
is the set of edges leaving $x$,
$$
\star(x)=\{ e\ |\ \iota(e)=x\}.
$$
\end{DEF}

The star of a vertex is the combinatorial analogue 
of the tangent space to a manifold at a point.

\begin{DEF}
A {\em connection} on a graph $\Gamma$ is a set of functions
$\nab_{(x,y)}$ or $\nab_e$, one for each oriented edge $e=(x,y)$ of
$\Gamma$, such that
\begin{enumerate}
\item $\nab_{(x,y)}: \star(x) \to \star(y)$,
\item $\nab_{(x,y)}(x,y) = (y,x)$, and
\item $\nab_{(y,x)} = (\nab_{(x,y)})^{-1}$.
\end{enumerate}
\end{DEF}
It follows that each $\nab_{(w,y)}$ is bijective, so each
connected component of $\Gamma$ is regular: all vertices have the same
valence.
Henceforth we will assume $\Gamma$ comes equipped with a connection
$\nab$.

\begin{DEF}
A {\em $3$-geodesic} is a sequence of four vertices $(x, y, z,w)$ with edges
$\{ x,y\}$, $\{ y,z\}$, and $\{ z,w\}$ for which
$\nab_{(y,z)}(y,x)=(z,w)$.  We inductively define a $k$-geodesic
as a sequence of $k+1$ vertices 
in the natural way. 
We may identify a
geodesic by specifying either its edges or its vertices, and we will
refer to {\em edge geodesics} or {\em vertex geodesics}
as appropriate.
The three
consecutive edges $(d,e,f)$ of a $3$-geodesic will be called an {\em
edge chain}.
\end{DEF}

\begin{DEF}
A {\em closed geodesic} is a sequence of edges $e_1,\dots ,e_n$ such that
each consecutive triple $(e_i,e_{i+1},e_{i+2})$ is an edge chain for
each $1\leq i\leq n$, modulo $n$.
\end{DEF}

A little care is required to understand when a geodesic is
closed, since it may in fact use some edges in  
$\star(x)$ multiple
times. It is not closed until
it returns to the same pair of edges in
the same order.
That is analogous to the fact that a periodic geodesic in
a manifold is an immersed submanifold, not an embedded
submanifold. The period completes only when it returns to a point with
the same velocity (tangent vector).

\begin{remark}
Because there is a unique closed geodesic through each pair of edges
in the star of a vertex, the set of all closed geodesics completely
determines the connection on $\Gamma$. 
We will sometimes use this fact
to describe a connection.
\end{remark}

We define totally geodesic subgraphs of a graph by analogy to
totally geodesic submanifolds of a manifold.

\begin{DEF}
Given a graph $\Gamma$ with a connection $\nab$, we say that a
subgraph $(V_0,E_0)=\Gamma_0\subseteq\Gamma$ is {\em totally geodesic}
if all geodesics starting in $E_0$ stay within $E_0$. 
\end{DEF}

\noindent This definition is equivalent to saying that a totally
geodesic subgraph $\Gamma_0$ is one in which, for every two adjacent
vertices $x$ and $y$ in $\Gamma_0$,  
$$
\nab_{(x,y)}(\star(x)\cap E_0)\subseteq E_0.
$$

Suppose now that $P =\{ e_1,\dots ,e_n\}$
is any cycle in $\Gamma$: $\tau(e_i)=\iota(e_{i+1})$ modulo
$n$. Then following the connection around
$P$ leads to a permutation
$$
\nab_{P}=\nab_{e_n}\circ\cdots\circ\nab_{e_1}
\circ\nab_{e_0} 
$$ 
of $\star(x)$.

\begin{DEF}
The {\em holonomy group} $Hol(\Gamma_x )$ at vertex $x$ of $\Gamma$ is
the subgroup of the permutation group of $\star(x)$ generated by
the permutations
$\nab_{P}$ for all cycles $P$ that pass through $x$.
\end{DEF}

It is easy to see that the holonomy groups $Hol(\Gamma_x )$ for the
vertices $x$ in each connected component of $\Gamma$ are
isomorphic. When $\Gamma$ is connected and $d$-regular we call that
group the holonomy group of $\Gamma$ and think of it as a subgroup of
$S_d$.

\section{Axial functions} \label{se:axial}

We described in Section~\ref{se:intro} how a graph arising from a GKM
manifold has associated to it an {\em axial function}, namely an
assignment of a vector in $\Z^n$ to each oriented edge $e$.  We need
to formalize this definition
for abstract graphs.

\begin{DEF}
An {\em axial function} on a graph with a connection is a map
$$
        \alpha : E \to \R^n\setminus\{ 0\}
$$
such that
$$
        \alpha(e^{-1}) = -\alpha(e)
$$
and for each $3$-geodesic $(d,e,f)$
$$
        \alpha(d) , \alpha(e), \alpha(f)
$$
are coplanar.
\end{DEF}

\noindent It follows immediately that the images under $\alpha$ of all
geodesics of $\Gamma$ are planar.
What matters about the axial function is the {\em direction} of
$\alpha(e)$ in $\R^n\setminus\{ 0\}$, not its actual value.  We
consider two axial functions $\alpha$ and $\alpha'$ to be equivalent
if
$$
\frac{\alpha(e)}{||\alpha(e)||}=\frac{\alpha'(e)}{||\alpha'(e)||}
$$ 
for all edges $e$.  Notice that $\alpha$ is {\em not} equivalent to
$-\alpha$.

If $e=(x,y)$ is an edge, we will denote $\alpha(e)$ by $\alpha(x,y)$,
rather than using two sets of parentheses.
We picture an edge chain as a succession of vectors joined head to
tail in their plane, as shown in the figure below.
\begin{figure}[h]
\begin{center}
\epsfig{figure=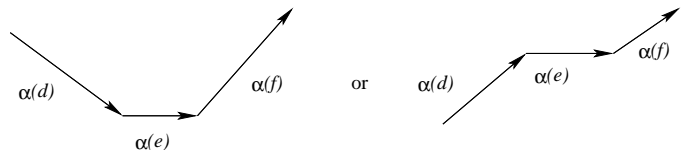,height=.75in}
\parbox{4.5in}{\caption[The GKM graph for CP3.]{{\small This shows how
we picture the axial function on an edge chain.}}}
\end{center}
\end{figure}
A picture of an equivalent axial function will show vectors with the
same orientations, but different lengths.

\begin{DEF}
An {\em immersion} of $(\Gamma,\alpha)$ is a map $F:V\to\R^n$ such
that
$$
\alpha(x,y)=F(y)-F(x).
$$
\end{DEF}

Our picture of an immersed vertex chain $(x,y,z,w)$ is shown below.
\begin{figure}[h]
\begin{center}
\epsfig{figure=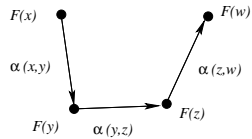,height=.75in}
\parbox{4.5in}{\caption[The GKM graph for CP3.]{{\small This shows how
we picture the axial function on an immersed vertex chain.}}}
\end{center}
\end{figure}
Here, the endpoints of the vectors do make sense, as the vertices are
points in $\R^n$.

\begin{DEF}
An axial function is {\em exact} if for every edge chain $(d,e,f)$,
\begin{equation}\label{eq:axialexact}
\alpha(f)+\alpha(d)= c\cdot \alpha(e),
\end{equation}
for $c\in\R$.
\end{DEF}

The figure below shows three exact edge chains.
\begin{figure}[h]
\begin{center}
\epsfig{figure=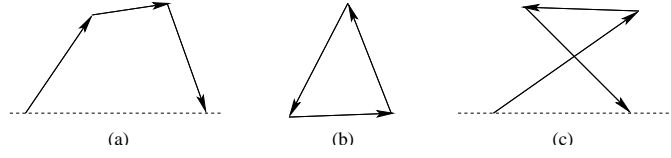,height=.75in}
\parbox{4.5in}{\caption[The GKM graph for CP3.]{{\small This shows
three exact edge chains.}}}
\end{center}\label{fig:exact}
\end{figure}
In Figure~4(a), the constant is positive, in (b) it is
zero, and in (c) it is negative.

When $\alpha$ is exact, the edges $\alpha(d)$ and $\alpha(f)$ must lie
on the same side of $\alpha(e)$ in the plan in which they lie.  The
examples in Section~\ref{subse:EGs} are immersible and exact.  So are
the axiam functions for the GKM manifolds discussed in
Section~\ref{subse:GKMgr}.  
Notice that the product of two
immersable (exact) axial functions is again immersable (exact).  Under
what conditions is an arbitrary axial function equivalent to one that
is immersable? exact?  We have only partial answers to these
interesting questions.

\begin{theorem}
If the axial function $\alpha$ is $3$-independent, then it
determines the connection.
\end{theorem}

\begin{proof}
Let $d$ and $e$ be edges with $\tau(d)=\iota(e)$.  Then
$3$-independence implies that there is only {\em one} edge $f$ with
$\tau(e)=\iota(f)$ 
and and $\alpha(f)$
in the plane determined by $\alpha(d)$ and $\alpha(e)$, 
\end{proof}

In Section~\ref{se:para}, we will need a special case of the following
definition.

\begin{DEF}
The product of two graphs $\Gamma_1=(V_1,E_1)$ and
$\Gamma_2=(V_2,E_2)$ is the graph
$$
\Gamma=\Gamma_1\times\Gamma_2=(V,E),
$$
with vertex set
$V=V_1\times V_2$. Two
vertices $(x_1,y_1)$ and $(x_2,y_2)$ are adjacent if and only if
\begin{enumerate}
\item $x_1=x_2$ and $\{ y_1,y_2\}\in E_2$; or
\item[(2)] $y_1=y_2$ and $\{ x_1,x_2\}\in E_1$.
\end{enumerate}
\end{DEF}

Suppose now that 
each 
$\Gamma_i$ is equipped with a connection $\nabla_i$
and axial function $\alpha_i:E_i\to\R^{n_i}$.
Then we can
define a connection $\nabla$ on $\Gamma$ in a natural way.  
by specifying the closed geodesics as the closed geodesics in each
component and some closed geodesics of
length $4$ which go between $\Gamma_1$ and $\Gamma_2$.  

The figure below shows one each of the two kinds of geodesics for
the example in which
$\Gamma_1$ is a $3$-cycle and $\Gamma_2$ is an
edge.  
\begin{figure}[h]
\begin{center}
\epsfig{figure=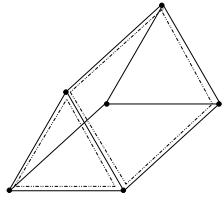,height=1in}

\parbox{4.5in}{\caption[The GKM graph for CP3.]{{\small This shows the
product of two graphs, 
showing one geodesic of each type.}}}
\end{center}
\end{figure}

We define an axial function
$\alpha:E\to\R^{n_1+n_2}$ by
$$
\alpha((x_1,y_1),(x_2,y_2))=(\alpha_1(x_1,x_2),\alpha_2(y_1,y_2)),
$$
where (by definition) $\alpha_i((x,x))=0$.  We leave it to the
reader to check that $\nabla$ is a well-defined connection, and that
$\alpha$ is indeed an axial function compatible with $\nabla$.  Note
that this generalizes the example of the hypercube, which is an
$n$-fold product of an edge.

\section{Betti numbers}\label{se:betti}

Suppose $\Gamma$ is a graph with a connection $\nabla$ and an axial
function $\alpha$ mapping edges to $\R^n$. The images under $\alpha$ of
the chains in $\Gamma$ are planar; we will study how those
chains wind in their planes.  To that end choose an arbitrary
orientation for each such plane $P$. Then whenever $\alpha(e) \in P$
the direction 
$\alpha(e)^{\perp}$ is a well defined direction in $P$. (If $\alpha$ is
immersible then 
$\alpha(e)^{\perp}$ is a well defined vector in $p$.)

Throughout this section we will assume $\alpha$ is
$2$-independent. That is, no two edges in the
star of a vertex of $\Gamma$ are mapped by $\alpha$ into the same line
in $\R^n$. Thus any two edges at a vertex determine a unique plane,
which we have assumed is oriented.

\begin{DEF}
The {\em curvature} $\kappa(d,e)$ of the edges
$d = (x,y)$ and $e = (y,z)$ at $y$ is
$\sign( \alpha(d)^{\perp}\cdot\alpha(e) ) \in\{ \pm 1\}$. 
\end{DEF}

\begin{DEF}
An edge chain $(d,e,f)$ is an {\em inflection} if $\kappa(d,e)$ and
$\kappa(e,f)$ have opposite signs. We say $\alpha$ is 
{\em inflection free} if there are no inflections.
\end{DEF}

Note that if an axial function is exact, then by \eqref{eq:axialexact}, it is
inflection free.

\begin{figure}[h]
\centerline{
\psfig{figure=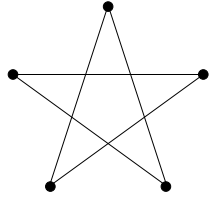,height=1in}
}
\centerline{
\parbox{4.5in}{\caption[Pentagram]{{\small The 
pentagram: a $5$-cycle with an unusual immersible axial function that is
nonetheless inflection free.}}}  }
\end{figure}

\begin{DEF}
A direction
$\xi\in \R^n\setminus\{ 0\}$ is {\em generic} if for all 
$e\in E$, 
$\xi\not\perp\alpha(e)$.
\end{DEF}

\begin{DEF}
The {\em index} of a
vertex $x\in V$ with respect to a generic direction $\xi$
is the number of edges $e\in\star(x)$ such that 
$$
\alpha(e)\cdot\xi <0.
$$ 
We call those the {\em down} edges.
Let $\beta_i(\xi)$ be the number of
vertices $x\in V$ such that the index of $x$ is exactly $i$.
\end{DEF}

\begin{theorem}
If $\Gamma$ is a graph with 
connection $\nab$ and an inflection free
axial function $\alpha$, then the Betti numbers $\beta_i$ do not
depend on the choice of direction $\xi$.
\end{theorem}

\begin{proof}
Imagine the direction $\xi$ varying continuously in $\R^n$. It is
clear from the definitions above that the indices of vertices can
change only when $\xi$ crosses one of the hyperplanes
$\alpha(x,y)^{\perp}$. Let us suppose that $(x,y)$ is the only edge of
$\Gamma$ at which the value of the axial function is a multiple of
$\alpha(x,y)$.  Then at such a crossing only the indices of the
vertices $x$ and $y$ can change. Suppose $\xi$ is near
$\alpha(x,y)^{\perp}$. Since $\alpha$ is inflection
free, the connection
mapping $\star(x)$ to $\star(y)$ preserves down edges, with the single
exception of edge $(x,y)$ itself. That edge is down for one of $x$ and
$y$ and up for the other. Thus the vertices $x$ and $y$ have indices
$i$ and $i+1$ for $\xi$ on one side of $\alpha(x,y)^{\perp}$ and
indices $i+1$ and $i$ on the other. Thus the number of vertices with
index $i$ does not change as $\xi$ crosses $\alpha(x,y)^{\perp}$.  If
there are several edges of $\Gamma$ at which the axial function is a
multiple of $\alpha(x,y)$, the same argument works, since by the
$2$-independence of $\alpha$, none of those edges can share a common
vertex.
\end{proof}

Henceforth we will assume $\alpha$ is inflection free.
The motivation for the following definitions comes from Morse theory.

\begin{DEF}
When the $\beta_i(\xi)$ are independent of $\xi$, we call them the 
{\em Betti numbers} of $\Gamma$ (or, more precisely, the Betti numbers
of the pair $(\Gamma, \alpha)$ ).
\end{DEF}

The following proposition is the combinatorial version of Poincar\'e
duality.

\begin{prop}\label{pr:poincare}
When the Betti numbers of a graph are independent of the choice of
$\xi$, then $\beta_i(\Gamma)=\beta_{d-i}(\Gamma)$ for $i=0,\dots,d$.
\end{prop}

\begin{proof}
Choose some $\xi$ with which to compute the Betti numbers of
$\Gamma$.  Then simply replace $\xi$ by $-\xi$, and a vertex of index
$i$ becomes a vertex of index $d-i$.
\end{proof}

\begin{DEF}
Given a generic $\xi$, a {\em Morse function compatible with $\xi$} on
a graph with an axial function $\alpha$ is a map $f:V\to\R$ such that if 
$(x,y)$ is an edge, $f(x)>f(y)$ whenever $\alpha(x,y)\cdot\xi>0$.
\end{DEF}

There is a simple necessary and sufficient condition for the existence
of a Morse function compatible with $\xi$.

\begin{theorem}
A Morse function compatible with $\xi$ exists if and only if there
exists no closed cycle $(e_1,\dots,e_n)$ with $e_1=e_n$, in $\Gamma$ for
which all the edges $e_i$ are ``up'' edges.
\end{theorem}

\begin{proof}
The necessity of this condition is obvious since $f$ has to be
strictly increasing along such a path.  To prove sufficiency, for
every vertex $p$, define $f(p)$ to be the length $N$ of the longest
path $(e_1,\dots,e_N)$ in $\Gamma$ 
of up edges
$\tau(e_N)=p$ 
The hypothesis that there is no cycle
of up edges guarantees that this function is well-defined, and it
is easy to check that it is a Morse function.
\end{proof}

\begin{remark}
One can easily arrange for $f$ in the above proof to be an {\em
injective} map of $V$ into $\R$ by perturbing it slightly.
\end{remark}

When $\alpha$ is immersible, so that $\alpha(x,y) = f(y) - f(x)$, then
we can 
define a Morse
function on $\Gamma$ by setting $m(x) = f(x) \cdot\xi$
for any generic direction $\xi$.
Then $m(x)$ increases
along each up edge.  The vertices with index $i$ resemble critical
points of Morse index $i$ in the Morse theory of a manifold. We call
the $\beta_i$ Betti numbers because when a graph we are studying
is the GKM graph of a manifold, the $\beta_i$ indeed correspond
to the Betti numbers of the manifold, and they are the dimensions of
the cohomology groups of the manifold.

\begin{remark}
When an inflection-free $2$-independent axial function is projected
generically into a plane, it retains those properties, so the Betti
numbers of $\Gamma$ can be computed using a generic direction in a
generic plane projection.  In most of our examples $\alpha$ is
immersible. In these cases we are of course drawing a planar embedding
of $\Gamma$. Thus the figures in this paper are more than mere
suggestions of some high dimensional truth. They actually capture all
the interesting information about $\Gamma$.
\end{remark}

\begin{DEF}
The generating function $\beta$ for the Betti numbers of $\Gamma$ is
the polynomial
$$
\beta(z)=\sum_{i=0}^{n}\beta_iz^i.
$$
\end{DEF}

\begin{remark}
When $\Gamma$ is $d$-regular, $\beta$ is of degree $d$. The sum of the
Betti numbers, $\beta(1)$, is just the number of vertices of $\Gamma$
\end{remark}

\begin{remark}
It is clear that $\beta_0 > 0$ if a Morse function exists, because the
vertex at which the Morse function assumes its minumum value has no
down edges.
\end{remark}

We can relate the Betti numbers of the product of two graphs to the
Betti numbers of the two multiplicands as follows.  The proof is left
to the reader.

\begin{prop}
Let $\Gamma$ and $\Delta$ be graphs with Betti numbers generated by
$\beta_\Gamma(z)$ and $\beta_\Delta(z)$ respectively.  Then
the generating function for the Betti numbers of the product
graph $\Gamma\times\Delta$ is the polynomial product
$$
\beta_\Gamma(z)\cdot\beta_\Delta(z).
$$
\end{prop}

\section{Simple convex polytopes}\label{se:simple}

Recall that a polytope $\Delta$ in $\R^n$ is simple if each vertex has
degree $n$. 
Every plane polygon is simple. Three of the
five Platonic 
polyhedra are simple: the tetrahedron, the cube and the
dodecahedron. So are the higher dimensional analogues of the simplex
and the cube we have already encountered.

Let $\Gamma$ be the embedded graph whose vertices and edges are
the one-skeleton of the simple polytope $\Delta$. Every two
dimensional face of $\Delta$ is a plane polygon. 
Since $\Delta$ is
simple, every pair of adjacent edges belongs to exactly one such
polygon. These polygons thus serve to define a natural connection on
$\Gamma$ whose geodesics are those polygons. That connection is
compatible with the exact axial function $\alpha$ defined by the
embedding. The totally geodesic subgraphs naturally correspond the
the faces of $\Delta$ of every dimension. If we assume $\Delta$ is
convex (and we shall) then $\alpha$ is inflection free and the Betti
numbers of $\Gamma$ are well defined.

Let us now turn to the identities (\ref{eq:1}). We will prove
that the $\beta_k$'s in these
identities are just the Betti numbers of $\Gamma$.  
Since the equations (\ref{eq:1}) can be solved for the $\beta_k$'s as
expressions in the terms on the left hand side, this will give us
another proof of Proposition~\ref{pr:poincare}.

\begin{proofof}{Proof of (\ref{eq:1})}
Fix a generic vector $\xi\in\R^n\setminus\{ 0\}$ ($\xi\not\perp\alpha(e)$
for all $e\in E$) and let $f:\R^n\to\R$ be the function $f(x)=x\cdot
\xi$.  By choosing $\xi$ appropriately, we can assume that $f:V\to\R$
is injective.  Then $f$ on $V$ is an
injective Morse function compatible with $\xi$ 
sense of Section~\ref{se:betti}.  Now to count the
faces of $\Delta$ we 
exploit the fact that every $n-k$-dimensional face $F$ has a unique
vertex $p_F$ at which $f$ takes its minimum value. 
At this vertex, the index of $p_F$ with respect to $\xi$ is at
most $k$, since the $n-k$ edges of $F$ at $p_F$ are upward-pointing
with respect to $\xi$.

Let $\ell$ be a number between $0$ and $k$, and let $p$ be a vertex of
$\Gamma$ of index $\ell$.  We can count the number of
$n-k$-dimensional faces $F$ for which $p_F=p$.  At $p$, there are exactly
$n-\ell$ edges which are upward pointing with respect to $\xi$, and
every $n-k$-element subset of this set of edges spans an
$n-k$-dimensional face $F$ with $p_F=p$.  Thus, there are
$\binom{n-\ell}{n-k}$ such faces in all.  Summing over $p$, we find
$$
\binom{n-\ell}{n-k}\beta_\ell
$$
$n-k$-dimensional faces $F$ for which the index of $p_F$
equals $\ell$.  Summing over $\ell$ we get
$$
f_{n-k} = \sum_{\ell=0}^k\binom{n-\ell}{n-k}\beta_\ell
$$
for the total number of $n-k$-dimensional faces of $\Delta$.
\end{proofof}

We noted in the introduction that Stanley proved the McMullen
conjectures by showing that the graph theoretical Betti numbers
$\beta_k$ in the identites (\ref{eq:1}) are actually the Betti numbers
of a complex projective variety, and hence that the identities
(\ref{eq:2}) follow from Poincar\'e duality and the
inequalities (\ref{eq:3}) from the hard Lefshetz
theorem.  We are part way there: the identites (\ref{eq:2}) are
Proposition~\ref{pr:poincare}.  A combinatorial proof of \eqref{eq:3}
exists (\cite{McMullen2}, \cite{Kalai}), but it is much harder than
the proof we gave above of \eqref{eq:2}.
So for completeness (and because it is so elegant) we 
will sketch Stanley's argument.

The McMullen inequalities follow from the following two assertions:
\begin{enumerate}
\item If $M$ is a GKM manifold, its manifold Betti numbers coincide
with its graph theoretical Betti numbers.

\item Given a simple convex polytope, there exists a GKM manifold
whose graph and axial function are the graph and axial function of the
polytope.
\end{enumerate}

The proof of the first of these assertions is by ``manifold'' Morse
theory and resembles our ``graph'' Morse theory computations of
Section~\ref{se:betti}. We will briefly describe how to prove the
second assertion. If we perturb the vertices of $\Delta$ slightly it
will still be simple and convex, so without loss of generality, we may
assume that the vertices of $\Delta$ are rational points in $\R^n$.
For each facet ({\em i.e.} $n-1$-dimensional face) $F$ of $\Delta$,
Then  we can find an outward pointing normal vector  $m_F$ to $\F$
with relatively prime integer entries. Let $F_1,\dots,F_d$ be the
facets of $\Delta$, and let $m_{F_i} = m_{1,i},\dots,m_{n,i})$.  The
map
$$
(\C^*)^d\to(\C^*)^n,
$$
mapping $(z_1,\dots,z_d)$ to $(w_1,\dots,w_n)$, where
$$
w_i=z_1^{m_{i,1}}\cdots z_d^{m_{i,d}},
$$
is a group homomorphism.  Let $N$ be its kernel  Clearly,
\begin{equation}\label{eq:quotient}
G=(C^*)^n=(\C^*)^d/N.
\end{equation}
Now let $(\C^*)^d$ act on $\C^d$ by coordinate-wise multiplication,
$$
a\cdot z=(a_1\cdot z_1,\dots,a_d\cdot z_d).
$$
Since $N$ is a subgroup of $(\C^*)^d$ it acts on
$\C^d$. Then (\ref{eq:quotient}) induces an action of $G$ on the
quotient
$$
\C^d/N.
$$
This quotient space, with its $G$ action, is, morally speaking, the
manifold $M$ that we are looking for.  Unfortunately, it is not a
manifold, and as a topological space, it is not even Hausdorff. 
Fortunately, it is easily desingularized, by deleting from it a finite
number of ``non-stable'' orbits of $G$.  More explicitly, for every
subset $I$ of $\{ 1,\dots,d\}$, let
$$
\C_I^d=\left\{ (z_1,\dots,z_d)\in\C^d\ |\ z_i=0\Leftrightarrow\ i\in
I\right\}.
$$
It is easy to see that these sets are exactly the $(\C^*)^d$-orbits in
$\C^d$.  In particular, if $F$ is a face of $\Delta$, there is a
unique subset $I$ of $\{ 1,\dots,d\}$ for which
$$
F=\cap_{i\in I}F_i
$$
and to that face, we attach the orbit
$$
\C^d_F:=\C^d_I.
$$
Finally, we define
\begin{equation}\label{eq:toricvar}
\C^d_{\Delta}=\cup_F\C^d_F.
\end{equation}

\begin{theorem}
The space $\C^d_{\Delta}$ is an open subset of $\C^d$ on which $N$
acts properly and locally freely, and the quotient space
$$
M=\C^d_{\Delta}/N
$$
is a compact orbifold.
\end{theorem}

\begin{remark}
The term ``orbifold'' means that $M$ is not quite a manifold.  It does
have singularities, but they are fairly benign, and its topological
properties, including the behavior of its Betti numbers, are the same
as those of a manifold.
\end{remark}

\begin{remark}
There is a simple condition equivalent to the assertion that
$M$ is a manifold: for
each vertex $p$ of $\Delta$ the $n$
vectors $m_{F_i}$ normal to the $n$ facets containing $p$
form a {\em
lattice} basis of $\Z^n$.  That is, every vector in $\Z^n$ can be
written as a linear combination of the $m_{F_i}$'s with integer coefficients.
\end{remark}

Since each of the summands in (\ref{eq:toricvar}) is a $(\C^*)^d$
orbit, the $G$-orbits in $M$ are the sets
$$
M_F=\C^d_F/N,
$$
and a simple computation shows that $\dim M_F=\dim F$.  Thus, the zero
and one dimensional orbits of $M$ correspond to the vertices and edges
of $\Delta$.  In other words, the action $\rho$ of $G$ on $M$ is a GKM
action and the graph of $(M,\rho)$ is the one-skeleton, $\Gamma$, of
$\Delta$.

\begin{EG}{Let $\Delta$ be the $n$-simplex
$$
\{ (x_1,\dots,x_n)\in\R^n_+\ |\ x_1+\cdots +x_n=1\}.
$$
Then
$$
m_{F_i}=(1,\dots,-n,\dots,1)
$$
(a ``$-n$'' in the $i^th$ slot), and $N$ is the diagonal subgroup of 
$(\C^*)^{n+1}$: the $n$-tuples for which $z_1=\cdots=z_{n+1}$.
Moreover, $\C_\Delta^{n+1}$
is $\C^{n+1}\setminus\{ 0\}$.  Note that the action of $N$ on
$\C^{n+1}$ is non-Hausdorff, as every $N$ orbit contains $0$ in its
closure.  However the action of $N$ on $\C^{n+1}\setminus\{ 0\}$ is
free and proper, and the quotient $\C_\Delta^{n+1}/N$ is just $\C
P^n$.  Thus, in this case, the construction we have just outlined
reproduces the first of the examples discussed in
Section~\ref{subse:EGs}.  For more information about this
construction, see, for example, \cite[p. 109--130]{G:toric}.}
\end{EG}

\section{Polynomial interpolation schemes}\label{se:poly}

\subsection{Polynomial interpolation schemes}
    \label{se:6.1}

Let $S$ be the polynomial ring in the variables
$x = (x_1,\ldots ,x_n)$ and $S^k$ the $k$\st{th}
graded component $S$: the space of homogeneous polynomials of
degree $k$.  
$S^k$ has dimension
$\binom{n+k-1}{n-1}$.

Let $\Gamma = (V,E)$ be a regular $d$-valent graph
with connection  $\triangledown$ and axial function
$\alpha: E \to  \R^n $. Then for every 
edge $e$, of $\Gamma$ we will
identify the vector, $\alpha(e) \in \R^n$, with the linear function
$\alpha_e(x) = \alpha(e) \cdot x$ so we can think of $\alpha_e$ as an
element of $S^1$. Finally, for $g$ and $h \in S$ we will say that
\begin{displaymath}
  g \equiv h \mod \alpha
\end{displaymath}
when $g-h$ vanishes on the hyperplane, $\alpha (x) =0$.

\begin{DEF}
Let $m$ be the numbers of vertices of $\Gamma$.  An $m$-tuple of polynomials 
\begin{displaymath}
  g_p \in S , \quad p \in V
\end{displaymath}
is a \emph{polynomial interpolation scheme} 
if for every $e = (p,q) \in E$ 
\begin{equation}
  \label{eq:6.1}
  g_p \equiv g_q \mod \alpha_e .
\end{equation}
Henceforth we will write $<g>$ for such an $m$-tuple.
We say that this scheme has \emph{degree} $k$ if for all $p$,
$g_p \in S^k$, and we will denote the set of all polynomial
interpolation schemes of degree $k$ by $H^k (\Gamma , \alpha)$.
\end{DEF}

It's clear that every interpolation scheme is in the space
\begin{displaymath}
  H^* (\Gamma ,\alpha) = \bigoplus^{\infty}_{k=0} H^k (\Gamma , \alpha)\, .
\end{displaymath}
Moreover this space is clearly a graded module over the ring $S$. That
is, if $<g_p>$ satisfies
(\ref{eq:6.1}) then for every $h \in S$, so does $<hg_p>$.
More generally, if $<g_p>$ and  $<h_p>$ satisfy
(\ref{eq:6.1}) then so does $<g_p\cdot h_p>$.
So $H^*(\Gamma ,\alpha)$ is not just a module, but in fact a
\emph{graded ring}, and $S$ sits in this ring as the subring of
constant interpolation schemes
\begin{displaymath}
  g_p =g  \hbox{  for all  } p.
\end{displaymath}

In this section we describe some methods for constructing
solutions of the interpolation equations \eqref{eq:6.1}.  These methods
will rely heavily on the ideas that we introduced in 
Sections~{\ref{se:conn} and \ref{se:axial}}.

Let $F: V \to \R^n$ be an immersion of $\Gamma$.
If we identify vector $F(p)$ with the monomial 
\begin{displaymath}
  f_p (x) = F(p) \cdot x
\end{displaymath}
then (\ref{eq:6.1}) is just a rephrasing of the identity
(\ref{eq:axialexact}), so $<f_p>$ is an
interpolation scheme of degree~1.  More generally if
\begin{displaymath}
  \p (z) = \sum^k_{i=0} \p_i(x) z^i
\end{displaymath}
is a polynomial in $z$ whose coefficients are polynomials in 
$x = (x_1 ,\ldots , x_n)$ then the $m$-tuple of polynomials
\begin{equation}
  \label{eq:6.2}
< \sum^k_{i=0} \p_i(x) f_p^i >
\end{equation}
is a polynomial interpolation scheme. $<f_p>$ itself corresponds to
the case $k=1$, $\p_0 = 0$, $p_1 = 1$.

\subsection{The complete graph}
\label{se:6.2}
In one important case this construction gives \emph{all} solutions of
(\ref{eq:6.1}).  Namely let $K_{n+1} = (V,E)$ be the complete graph on
$n+1$ vertices with the natural connection. Then every immersion
$F: V \to \R^n$ defines an axial function compatible with that
connection by setting
\begin{equation}
  \label{eq:6.3}
  \alpha(p,q) = F(q) - F(p)
\end{equation}
for every oriented edge $e = (p,q)$.  We will prove

\begin{theorem}\label{th:thof6.2}
  If the axial function (\ref{eq:6.3}) is two-independent, every
interpolation scheme can be written \emph{uniquely} in the form
\begin{equation}
  \label{eq:6.4}
	<g_p> = < \sum^k_{i=0} \p_i(x) f_p^i >
\end{equation}
for some polynomials $\p_i$.
\end{theorem}

\begin{proof}
  By induction on $n$.  Let $\{ p_1 , \ldots , p_{n+1} \}$ be the
vertices of $V$, and let $<g_p>$ be an
interpolation 
scheme.  By induction there exists a polynomial
  \begin{displaymath}
    \p (z) = \sum^k_{i=0} \p_i z^i
  \end{displaymath}
with coefficients in $S$ such that $\p (g_{p_i}) = f_{p_i}$ for $i=1,
\ldots , n$.  Hence the $n+1$-tuple of polynomials 
\begin{displaymath}
  <f_p -p(g_p)>n
\end{displaymath}
is an interpolation scheme vanishing on $p_1 , \ldots , p_{n}$.
Therefore, by (\ref{eq:6.1}) and the two-independence of the axial
function (\ref{eq:6.2}) 
  \begin{displaymath}
    f_{p_n} -\p (g_{p_n}) = h \prod_{i<n} (g_{p_n} - g_{p_i})
  \end{displaymath}
for some polynomial $h \in S$.  Let
\begin{displaymath}
  \q(z) = h \prod_{i \le n} (z-g_{p_i}) \, .
\end{displaymath}
Then
\begin{eqnarray*}
  \q (g_{p_{n+1}}) &=& f_{p_{n+1}} -\p (g_{p_{n+1}})\\
\noalign{\hbox{and}}
\q (g_{p_{n+1}}) &=& 0 
\end{eqnarray*}
for $i<n$.  
Thus the theorem is true for $\# V=n+1$ with $\p$ replaced by $\p
+\q$.

The uniqueness of $\p$ follows from
the Vandermonde identity 
\begin{eqnarray*}
  \det \left(
    \begin{array}{cccc}
      1 & g_{p_1} & \ldots & g^{n-1}_{p_1}\\
      \vdots & \vdots & \ddots & \vdots \\
      1&  g_{p_n}& \ldots & g^{n-1}_{p_n}
    \end{array}
\right) = \prod_{i>j} g_{p_i} - g_{p_j} \, ,
\end{eqnarray*}
the right-hand side of which is non-zero by the two-independence of the axial function (\ref{eq:6.2}).

\end{proof}

\subsection{Holonomy and polynomial interpolation schemes}
\label{se:6.3}

The complete graph is the only example we know of for which the
methods of the previous section
give \emph{all} the polynomial interpolation schemes.
In this section we describe an 
alternative method which is effective in examples in which one
has information about the holonomy group of the graph $\Gamma$.
To simplify the exposition below we will confine ourselves to the case
in which the axial function $\alpha$ is exact.

Let $p_0$ be a vertex of $\Gamma$.  The holonomy group, 
$\Hol(\Gamma_{p_0})$ is by definition a subgroup of the group of
permutations of the elements of $\star (p_0)$, so if we enumerate its
elements in some order 
\begin{displaymath}
  e^0_i \in \star (p_0) \quad i=1, \ldots d
\end{displaymath}
we can regard $\Hol (\Gamma_0)$ as a subgroup of the permutation group
$S_d$ on $\{1, \ldots, d\}$.
Let $\q (z_1 , \ldots , z_d)$ be a polynomial
in $d$ variables with scalar coefficients.  We will say that $\q$ is
$\Hol (\Gamma_{p_0})$ \emph{invariant} if for every $\sigma \in \Hol
(\Gamma_{p_0})$
\begin{displaymath}
  \q (z_{\sigma(1)} , \ldots , z_{\sigma (n)}) = \q (z_1 , \ldots , z_n)\, .
\end{displaymath}
Now fix such a $\q$ and construct a polynomial assignment $<g_p>$
as follows.  Given a path, $\gamma$ in $\Gamma$ joining $p_0$ to $p$
the connection gives us a holonomy map 
\begin{displaymath}
  \triangledown_{\gamma} : \star (p_0) \to \star (p)
\end{displaymath}
mapping $e^0_1, \ldots , e^0_d$ to $e_1, \ldots , e_d$.  Set
\begin{equation}
  \label{eq:6.5}
  g_p = \q (\alpha_{e_1} (x) , \ldots , \alpha_{e_d} (x)) \, .
\end{equation}
The invariance of $\q$
guarantees that this definition is independent of the choice of
$\gamma$.  Let us show that $<g_p>$ satisfies the
interpolation conditions (\ref{eq:6.1}).  Let $e=(p,q)$.
The map
$\Gamma$
 \begin{displaymath}
   \triangledown_p : \star (p) \to \star (q)
 \end{displaymath}
maps $e_1 , \ldots , e_d$ to  $e'_1 , \ldots , e'_d$  and by the exactness of $\alpha$
\begin{displaymath}
  \alpha_{e'_i} \equiv \alpha_{e_i} \mod \alpha_e \, .
\end{displaymath}
Hence
\begin{displaymath}
  \q (\alpha_{e'_1}, \ldots ,\alpha_{e'_d}) \equiv
  \q (\alpha_{e_1}, \ldots ,\alpha_{e_d}) \mod \alpha_e \, .
\end{displaymath}

If the holonomy group is small this construction provides
many solutions of (\ref{eq:6.1}).
Even if $\Hol (\Gamma_{p_0})$
is large this method yields some interesting solutions.
For instance if $\q$ is a symmetric polynomial in
$z_1, \ldots , z_d$, (\ref{eq:6.5}) is a solution of  (\ref{eq:6.1}).

\subsection{Totally geodesic subgraphs and polynomial interpolation schemes}
\label{se:6.4} 
A third method for constructing solutions of (\ref{eq:6.1}) makes use of
totally geodesic subgraphs.
Whenever $\Gamma_0 =
(V_0,E_0)$ is a totally geodesic subgraph of degree $j$
then for every $p \in V_0$, 
$\star (p)$ is a disjoint union of $\star (p,\Gamma_0)$
and its complement, which we can
regard as the \emph{tangent} and \emph{normal}
spaces to $\Gamma_0$ at $p$.
Let
\begin{equation}
  \label{eq:6.6}
  g_p = \prod_{e \perp \Gamma_0} \alpha(e) \ ,
\end{equation}
a homogeneous polynomial of degree $d-j$.
By the results of Section~\ref{se:6.3}, 
the assignment $<g_p>$ sending $p \to g_p$, is 
an interpolation scheme on $V_0$, and we can extend this scheme to $V$
by setting
\begin{equation}
  \label{eq:6.7}
  g_p = 0 
\end{equation}
for $p \in V-V_0$.  Then (\ref{eq:6.6}) and (\ref{eq:6.7})
\emph{do} define an interpolation scheme on $V$. Clarly the
interpolation conditions (\ref{eq:6.1}) are satisfied if $e=(p,q)$ is
either 
an edge of $V_0$ or if $p$ and $q$ are both in
$V-V_0$.  If $p \in V_0$ and $q \in V-V_0$
then $\alpha(p,q)$ is one of the factors in the product
(\ref{eq:6.6}); so in this case the interpolation condition (\ref{eq:6.1})
is also satisfied.

\section{Polynomial interpolation schemes and Betti numbers}\label{se:para}

In many examples, all solutions of the interpolation equations
\eqref{eq:6.1} can be constructed by the methods discussed in
Section~\ref{se:poly}. However, to check this, one needs an effective
way of counting the total number of solutions of these equations.
Modulo some hypotheses which we will make explicit below, we will
prove the following theorem.

\begin{theorem}\label{th:interp}
The dimension of the space of solutions, $H^r(\Gamma,\alpha)$, of the
interpolation equations \eqref{eq:6.1} is equal to
\begin{equation}\label{eq:7.1}
\sum_{\ell=0}^r \binom{r-\ell+n-1}{n-1}\beta_\ell.
\end{equation}
\end{theorem}

We will show in the next section that
the resemblance between \eqref{eq:7.1}
and the McMullen formula \eqref{eq:1}
is not entirely an accident.

Let $\xi\in\R^n$ be a generic vector.  Without loss of generality, we may
assume that $\xi$ is the vector $(1,0,\dots,0)$.  Hence,
\begin{equation}\label{10.1}
\alpha(e)=m_e\cdot(x_1-\delta(e)),
\end{equation}
where $m_e\neq 0$ and 
\begin{equation}\label{10.2}
\delta(e)=d_{e,2}x_2+\cdots+d_{e,n}x_n.
\end{equation}
Our proof of Theorem~\ref{th:interp} will make three unnecessarily
stringent assumptions about $\Gamma$, $\nabla$, and $\alpha$.  That
these assumptions {\em are} unnecessarily stringent is shown in
\cite{GZ2}.  The proof given in \cite{GZ2} of Theorem~\ref{th:interp}
follows the lines of the proof below, but is much more complicated.

The three assumptions that we will make are the following.  The first
is that the axial function $\alpha$ is $3$-independent.  That is, for
each $p\in V$, the set of vectors
\begin{equation}\label{10.3}
\alpha(\star(v))
\end{equation}
is $3$-independent: any set of three of these vectors is linearly
independent.  This assumption guarantees that every plane in $\R^n$
will contain the image under $\alpha$ of a disjoint set of closed
geodesics.
Namely, let $W$ be a two-dimensional subspace of $\R^n$ and let
$\Gamma_W=(V,E_W)$ be the subgraph whose oriented edges $e\in E_W$
satisfy $\alpha(e)\in W$.  Then $3$-independence guarantees that if
$\Gamma_W$ is non-empty, its
connected components are either edges or  totally geodesic subgraphs of
valence two.  In particular, if $e_1$ and $e_2$ are edges of $\Gamma$
having a common vertex $p$, there is a {\em unique} connected totally
geodesic subraph of valence two containing $e_1$ and $e_2$.  To see
this, let $W=\spann\{\alpha(e_1),\alpha(e_2)\}$ and let $F$ be the
connected component of $\Gamma_W$ containing $p$. 

To avoid repeating the phrase ``connected totally geodesic subgraphs
of valence two,'' we will henceforth refer to such objects as {\em
two-faces}.  Our second assumption about $\Gamma$, $\nabla$ and
$\alpha$ is that for every two-face $F$, the zeroth Betti number
$\beta_0(F)$ is one.  In other words, if a totally geodesic subgraph
of valence two is connected in the graph theoretical sense, it is
connected in the homological sense as well.

Finally we will asume that there exists a $\xi$-compatible Morse
function,
\begin{equation}\label{10.4}
f:V\to\R,
\end{equation}
as defined in Section~\ref{se:betti}. By perturbing $f$ slightly, we
may assume that the map (\ref{10.4}) is injective.

\begin{DEF}
The {\em critical values} of $f$ are the numbers $f(p)\in\R$ for $p\in
V$.  The {\em regular values} of $f$ are the numbers which are {\em
not} critical values.
\end{DEF}

Our goal for the remainder of this section is to prove the following
theorem.

\begin{theorem}\label{th:||redrawing}
Suppose $\Gamma$ is a graph equipped with a connection $\nabla$, a
3-independent axial function $\alpha:E\to\R^n$, and a $\xi$-compatible
Morse function $f$.  Suppose further that each two-face of $\Gamma$
has zeroth Betti number equal to $1$.  Then the dimension of $H^r(\Gamma,
\alpha)$ is given by the formula \eqref{eq:7.1}.
\end{theorem}

\subsection{The cross sections of a graph}
Let $c\in \R$ be a regular value of $f$, and let $V_c$ be the set of
oriented edges $e$ of $\Gamma$ having $f(\iota(e))<c<f(\tau(e))$.
Intuitively, these are the edges of $\Gamma$ which intersect the
``hyperplane'' $f=c$.  The main goal of this section is the
construction of the {\em cross section of $\Gamma$ at $c$} which will
be a graph $\Gamma_c$ having $V_c$ as its vertex set.  This
construction will make strong use of the second of our three hypotheses.
This hypothesis implies the following.

\begin{lemma}
For every two-face, $F$, the restriction of $f$ to $V_F$ has a unique
maximum, $p$, and a unique minimum, $q$.  Moreover, if $f(p)>c$ and
$f(q)<c$, there are exactly two edges of $F$ contained in $V_c$.
\end{lemma}

We will now decree that a pair $(e_1,e_2)\in V_c\times V_c$ is in $E_c$
if and only if they are contained in a common two-face.  The following
is an immediate corollary of the lemma.

\begin{theorem}
The pair $\Gamma_c=(V_c,E_c)$ is a $(d-1)$-regular graph.
\end{theorem}

\begin{remark}
One can intuitively think of the edges of $\Gamma_c$ as being the
intersections of two-faces with the ``hyperplane'' $f=c$.  In this
way, each vertex of $\Gamma_c$ corresponds to an edge of $\Gamma$, and
each edge of $\Gamma_c$ corresponds to a two-face of $\Gamma$.
\end{remark}

We will now show how to equip $\Gamma_c$ with a connection and an
axial function.  Let $F$ be a two-face, $p$ a vertex of $F$ and $e_1$
and $e_2$ the two edges of $F$ meeting in $p$.  We define two natural
connections on $\Gamma_c$.

The first connection is the {\em up connection} $\nabla_{up}$ on
$\Gamma_c$.  Let $\tau(e_1)=q_1$ and $\tau(e_2)=q_2$.  Suppose $e\neq
e_1$ is an edge of $\Gamma$ with $\iota(e)=\tau(e_1)$.  If $q$ is the
maximum point of $f$ on $F$, then there is a unique geodesic path on
$F$ joining $q_1$ to $q_2$ passing through $q$.  By applying $\nabla$
to $e$ along this path, we can associate to $e$ an edge $e'$ of
$\Gamma$ with $\tau(e')=q_2$.  Let $E$ be the unique two-face
containing $e_1$ and $e$, and $E'$ the unique two-face containing
$e_2$ and $e'$.  Then the correspondence $\nabla_{up}$ mapping
$E\mapsto E'$ defines a connection on $\Gamma_c$.

We next define the {\em down connection} $\nabla_{down}$ on $\Gamma_c$.  Let
$\iota(e_1)=p_1$ and $\iota(e_2)=p_2$.  Suppose $e\neq e_1$ is an edge
of $\Gamma$ with $\iota(e)=\iota(e_1)$.  If $p$ is the minimum point of
$f$ on $F$, then there is a unique geodesic path on $F$ joining $p_1$
to $p_2$ passing through $p$.  By applying $\nabla$ to $e$ along this
path, we can associate to $e$ an edge $e'$ of $\Gamma$ with
$\iota(e')=p_2$.  Let $E$ be the unique two-face containing $e_1$ and
$e$, and $E'$ the unique two-face containing $e_2$ and $e'$.  Then the
correspondence $\nabla_{down}$ mapping $E\mapsto E'$ defines a
connection on $\Gamma_c$.

Now we will give a natural axial function on $\Gamma_c$.  Let $W$ be
the two-dimensional space $\spann\{\alpha(e_1),\alpha(e_2)\}$.

\begin{lemma}\label{le:10.2}
For every edge $e$ of $F$, $\alpha(e)\in W$.
\end{lemma}

\begin{proof}
Since $F$ is totally geodesic, it is a connected component of
$\Gamma_W$.  Furthermore, since $F$ has valence $2$, it must, in fact,
be a geodesic, and so its image with respect to $\alpha$ must lie in a
plane.  This completes the proof.
\end{proof}

We will identify $\R^{n-1}$ with the orthogonal complement of
$\xi=(1,0,\dots,0)$ in $\R^n$ and let $\delta_F$ be a fixed basis
vector of the one-dimensional space $W\cap \R^{n-1}$.  For every edge
$e$ of $\Gamma$, let $\delta(e)$ be the vector in $\R^{n-1}$ defined by
(\ref{10.1}) and (\ref{10.2}).

\begin{lemma}\label{le:10.3}
For every pair of edges $e$ and $e'$ of $F$, $\delta(e)-\delta(e')$
is a multiple of $\delta_F$.  
\end{lemma}

\begin{proof}
By Lemma~\ref{le:10.2}, $\alpha(e)$ and $\alpha(e')$ are in $W$, so
\begin{equation}\label{eq:*}
\delta(e')-\delta(e)=m_{e'}^{-1}\cdot \alpha(e')-m_e^{-1}\cdot
\alpha(e)
\end{equation}
is in $W\cap\R^{n-1}$.
\end{proof}

\begin{theorem}
Let $e_1$ and $e_2$ be adjacent elements of $V_c$ and let $F$ be the
unique two-face with $e_1$ and $e_2$ as edges. The assignment
$(e_1,e_2)\mapsto \delta((e_1,e_2))=\delta_F$ is an axial function on
$\Gamma_c$ taking its values in $\R^{n-1}$.  This is an axial function
compatible with $\nabla_{up}$ {\em and} with $\nabla_{down}$.
\end{theorem}

\begin{proof}
We will prove that this is an axial function compatible with
$\nabla_{up}$.  The proof that it is compatible with $\nabla_{down}$
is nearly identical.

Let $\tau(e_1)=q_1$ and $\tau(e_2)=q_2$ and let $q$ be the unique
maximum point of $f$ on $F$.  Suppose $e\neq e_1$ is an edge of
$\Gamma$ with $\iota(e)=\tau(e_1)$.  In this case, we apply $\nabla$
to $e$ along the unique path from $q_1$ to $q_2$ through $q$, and
associate to $e$ an edge $e'$ of $\Gamma$ with $\tau(e')=q_2$.  If $E$
is the unique two-face containing $e_1$ and $e$, and $E'$ is the
unique two-face containing $e_2$ and $e'$, then $\nabla_{up}$ maps $E$
to $E'$.  Thus, $(E,F,E')$ is an edge chain under $\nabla_{up}$ in
$\Gamma_c$.  The three dimensional subspace of $\R^n$ spanned by
$\alpha(e_1)$, $\alpha(e_2)$ and $\alpha(e)$ is the same as the three
dimensional subspace of $\R^n$ spanned by $\alpha(e_1)$, $\alpha(e_2)$
and $\alpha(e')$, since $e$ was obtained from $e'$ by the original
connection $\nabla$.  Hence the intersections of these
three-dimensional subspaces with $\R^{n-1}$ are the same
two-dimensional subspaces.  Therefore, the image of $(E,F,E')$ under
$\delta$ lies in a plane, and thus $\delta$ is an axial function
compatible with $\nabla_{up}$.
\end{proof}

In Section~\ref{subse:morseineq}, we will need the following lemma.

\begin{lemma}\label{le:10.4}
The axial function $\delta$ is 2-independent.
\end{lemma}

\begin{proof}
Let $e$ be in $V_c$ and let $e_1$ and $e_2$ be edges of $\Gamma$
with $p=\tau(e)=\iota(e_i)$.  Let $F_i$ be the two-face with edges $e$
and $e_i$, for $i=1,2$.  Then up to scalar multiple,
$$
\beta_{F_i}=\beta(e_i)-\beta(e)
$$
by Lemma~\ref{le:10.3}, and by 3-independence, $\alpha(e)$,
$\alpha(e_1)$, and $\alpha(e_2)$ are linearly independent.  Hence, by
(\ref{10.1}) and by (\ref{10.2}), $\beta_{F_1}$ and $\beta_{F_2}$ are
linearly independent.
\end{proof}

\subsection{A Morse lemma for cross sections.}\label{subse:7.2}

A classical theorem in Morse theory describes how the level sets of a
Morse function change as one passes through a critical point.  The
goal of this section is to prove a combinatorial analogue of this
theorem.

\begin{theorem}\label{th:MorseLemma1}
Let $p\in V$ be a vertex of index $k$ and let $c=f(p)$ and
$c^{\pm}=c\pm\varepsilon$.  Then for small $\varepsilon$,
$\Gamma_{c^+}$ can be obtained from $\Gamma_{c^-}$ by deleting a complete
totally geodesic subgraph isomorphic to $K_k$ and inserting in its
place a complete totally geodesic subgraph isomorphic to $K_{d-k}$.
\end{theorem}

\begin{remark}
In the previous section, we defined {\em two} canonical connections on
a cross section, the up connection and the down connection.  The
subgraph that we delete from $\Gamma_{c^-}$ will be totally geodesic with
respect to the up connection, and the subgraph of $\Gamma_{c^+}$ that
we insert in its place will be totally geodesic with respect to the
down connection.
\end{remark}

\begin{proof}
We will orient the edges of $\Gamma$ by assigning to each edge the
orientation for which $\alpha(e)\cdot\xi>0$.  Thus, because $f$ is
compatible with our choice of $\xi$,
$$
f(\tau(e))>f(\iota(e)).
$$
Let $e_1,\dots,e_k$ be the oriented edges with $\tau(e_r)=p$ and
$e_1',\dots,e_{\ell}'$ be the oriented edges with $\iota(e_s')=p$.
Here, $\ell=d-k$, and $k$ is by definition the idex of $p$.  Then, if
$\varepsilon$ is sufficiently small, $f(\iota(e_r))<c^-$ and
$f(\tau(e_s'))>c^+$.  See the figure below.
\begin{figure}[h]
\centerline{
\epsfig{figure=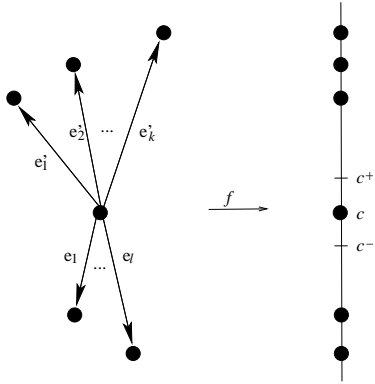,height=2in}
}
\centerline{
\parbox{4.5in}{\caption[Morse lemma]{{\small This
shows the edges appearing as vertices in $\Gamma_{c^-}$ and in
$\Gamma_{c^+}$.}}}
}
\end{figure}
Thus, the edges $e_1,\dots,e_k$ are vertices of $\Gamma_{c^-}$ and the
edges $e_1',\dots,e_{\ell}'$ are vertices of $\Gamma_{c^+}$.  By
3-independence, there exists a unique two-face $F_{i,j}$ having $e_i$
and $e_j$ as edges, so the $e_i$ regarded as vertices of
$\Gamma_{c^-}$ form a subgraph $\Delta_-$ of $\Gamma_{c^-}$ isomorphic
to $K_k$.  Moreover, the ``up'' geodesic path in $F_{i,j}$ joining
$\tau(e_i)=p$ to $\tau(e_j)=p$ consists of the point $p$ itself.  So,
$\nabla$ on $\Gamma$ simply maps $e_m$, $m\neq i,j$, along the path
$p$ to itself, and thus 	$\nabla_{up}$ maps the two-face $F_{i,m}$ to
the two-face $F_{j,m}$.  Thus, $\Delta_-$ is a totally geodesic subgraph
of $\Gamma_{c^-}$ with respect to the up connection.  Note that this
connection on $K_k$ agrees with the usual connection on $K_k$, as
described in Section~\ref{se:EGs}.

Similarly, the edges, $e_1',\dots,e{\ell}'$ are the vertices of a
subgraph $\Delta_+\iso K_{\ell}$ of $\Gamma_{c^+}$ which is totally
geodesic with respect to the down connection.

Now let $e$ be an edge of $\Gamma$ not having $p$ either as an initial
or terminal vertex and satisfying 
$$
f(\iota(e))<c<f(\tau(e)).
$$
If we choose $\varepsilon$ small enough, all of these edges $e$ belong
to both $V_{c^-}$ and $V_{c^+}$.  Thus, there is a bijection
$$
V_{c^-}-\{e_i\ |\ i=1,\dots,k\}\to V_{c^+}-\{e_i'\ |\
i=1,\dots,\ell\}.
$$
In other words, $\Gamma_{c^+}$ is obtained from $\Gamma_{c^-}$ by
deleting $\Delta_-$ and inserting $\Delta_+$.
\end{proof}

\subsection{A Morse lemma for $H^r(\Gamma_c)$}\label{subse:7.3}

We will compute in this section the change in dimension of
$H^r(\Gamma_c)$ as $c$ passes through a critical value of indes $k$.
As in Subsection~\ref{subse:7.2}, let $c=f(p)$ and let
$\Gamma_{c^\pm}=\Gamma_\pm$ be the cross-sections of $\Gamma$ just
above and just below $f=c$.  Let $V_\pm$ be the vertices of
$\Gamma_\pm$ and $V_\pm^c$ the vertices of the subgraphs $\Delta_\pm$
of $\Gamma_\pm$.  We will prove the following change of dimension
formula.

\begin{lemma}
\begin{equation}\label{eq:7A}
\dim H^r(\Gamma_+)-\dim H^r(\Delta_+)=\dim H^r(\Gamma_-)-\dim
H^r(\Delta_-).
\end{equation}
\end{lemma}

\begin{proof}
An element of $H^r(\Gamma_+)$ is a function which assigns to each
vertex $e$ of $\Gamma_+$ a homogeneous polynomial in $x_2,\dots,x_n$
of degree $r$ and for each edge of $\Gamma_+$ satisfies the $\Gamma_+$
analogue of the interpolation conditions \eqref{eq:6.1}.  Therefore,
since $\Delta_+$ is a totally geodesic subgraph of $\Gamma_+$, the
restriction of this function to $V_+^c$ is an element of
$H^r(\Delta_+)$, so there is a natural mapping
\begin{equation}\label{eq:7B}
H^r(\Gamma_+)\to H^r(\Delta_+).
\end{equation}
We claim that this mapping is surjective.  Indeed by
Theorem~\ref{th:thof6.2} and the formula \eqref{eq:*}, every element
of $H(\Delta_+)$ is the restriction to $V_+^c$ of an element of
$H(\Gamma_+)$ of the form $\p(\delta)$ where $\delta$ is given by
\eqref{10.2} and 
$$
\p=\p(z)=\sum_{i=1}^{\ell-1}\p_iz^i
$$
is a polynomial in $z$ of degree $\ell-1$, where $\ell=d-k$, having as
coefficients polynomials $\p_i$ in $x_2,\dots,x_n$.

The analogous mapping 
\begin{equation}\label{eq:7C}
H^r(\Gamma_-)\to H^r(\Delta_-).
\end{equation}
is also surjective, so if we denote the kernels of \eqref{eq:7B} and
\eqref{eq:7C} by $H^r(\Gamma_+,\Delta_+)$ and $H^r(\Gamma_-,\Delta_-)$
respectively, the proof of \eqref{eq:7A} reduces to showing that
\begin{equation}\label{eq:7D}
\dim H^r(\Gamma_+,\Delta_+)=\dim H^r(\Gamma_-,\Delta_-).
\end{equation}
By Theorem~\ref{th:MorseLemma1},
$$
V_+-V_+^c=V_--V_-^c,
$$
so if
\begin{equation}\label{eq:7E}
f_p^-, \ p\in V_-
\end{equation}
is an interpolation scheme that vanishes on $V_-^c$ we can associate
with it a function
\begin{equation}\label{eq:7F}
f_p^+,\ p\in V_+
\end{equation}
by setting
\begin{equation}\label{eq:7G}
f_p^+=f_p^-
\end{equation}
for $p\in V_=-V_+^c$ and
\begin{equation}\label{eq:7H}
f_p^+=0
\end{equation}
for $p\in V_+^c$.  Let us show that \eqref{eq:7F} is an interpolation
scheme for the graph $\Gamma_+$.  To prove this, we must show that 
$$
f^+_{e_1}=f_{e_2}^+ \mod \beta_F
$$
for every edge $F=(e_1,e_2)$ of $\Gamma_+$.  If $e_1$ and $e_2$ are in
$V_+-V_+^c$, this follows from \eqref{eq:7G}.  If $e_1$ and $e+2$ are
in $V_+^c$, then if follows from \eqref{eq:7H}.  So the only case we
have to consider is the case where $e_1\in V_+^c$ and $e_2\in
V_+-V_+^c$.  For the moment, let us regard $e_1$ and $e_2$ as edges of
$\Gamma$ and $F$ as a $2$-face of $\Gamma$ containing these edges.
The critical point $p$ on the level set $f=c$ is a vertex of $F$ since
$p$ is the initial vertex of $e_1$ and the two edges of $F$ meeting in
$p$ must be $e_1$ and one of the $e_j$'s in $V_-^c$, since if it were
in $V_+^c$, $e_2$ would have to be that $e_j$.  Now note that
$$
f_{e_2}^+=f_{e_2}^-
$$
by \eqref{eq:7G} and
$$
f_{e_1}^+=f_{e_j}^-=0
$$
by \eqref{eq:7H} and by the fact that $f^-$ is zero on $V_-^c$.
Hence,
$$
f_{e_2}^+-f_{e_1}^+=f_{e_2}^--f_{e_j}^- =0\mod\delta_F,
$$
since $e_2$ and $e_j$ lie on the common edge $F$ of $\Gamma_-$.  Thus,
the natural map \eqref{eq:7B} is indeed a surjection.

The dimensions of $H^r(\Delta_+)$ and $H^r(\Delta_-)$ can be computed
directly from Theorem~\ref{th:MorseLemma1}.  Namely, by
Theorem~\ref{th:MorseLemma1}, the dimension of $H^r(\Delta_+)$ is the
dimension of the space of homogeneous polynomials of degree $r$ in
$z,x_2,\dots,x_n$ of the form
$$
\sum_{i=0}^{\ell-1} \p_i(x_2,\dots,x_n)z^i,\ \ell=d-k,
$$
so if we let $d_r(n)$ be the dimension of the space of homogeneous
polynomials in $x_1,\dots,x_n$ of degree $r$, we get the formula
\begin{equation}\label{eq:7J}
\dim H^r(\Delta_+)=\sum_{i=1}^{\ell-1} d_{r-i}(n-1).
\end{equation}
Similarly, we get the formula
\begin{equation}\label{eq:7K}
\dim H^r(\Delta_-)=\sum_{i=1}^{k-1} d_{r-i}(n-1).
\end{equation}
Noting that
$$
d_s(n)=\sum_{i=0}^sd_i(n-1),
$$
we can rewrite \eqref{eq:7J} and \eqref{eq:7K} as
\begin{equation}\label{eq:7J'}
\dim H^r(\Delta_+)= d_{r}(n)-d_{r-\ell}(n)
\end{equation}
and
\begin{equation}\label{eq:7K'}
\dim H^r(\Delta_-)= d_{r}(n)-d_{r-k}(n).
\end{equation}
Hence, from \eqref{eq:7A}, we get the identity
\begin{equation}\label{eq:7L}
\dim H^r(\Gamma_+)-\dim H^r(\Gamma_-)=d_{r-k}(n)-d_{r-\ell}(n).
\end{equation}
\end{proof}

\subsection{The proof of Theorem~\ref{th:||redrawing}.}

We are now ready to prove the main theorem of this section.

\begin{proofof}{Proof of Theorem~\ref{th:||redrawing}}
Regard the unit interval $I=[0,1]$ as the graph consisting of a single
edge.  This graph has a unique connection.  Furthermore, equip it with
the axial function $\alpha(0)=x$ and $\alpha(1)=-x$, where $x$ is the
unit vector $1$ in $\R$.  Now let $\Gamma$ be a graph with a
connection and axial function, and let $\tilde{\Gamma}$ be the product
graph $\Gamma\times I$ with its product axial function
$\tilde{\alpha}: \tilde{\Gamma}\to\R^n\times\R=\R^n$, as defined in
Section~\ref{se:axial}.  Let $f: V_{\Gamma}\to\R$ be our given Morse
function, and extend $f$ to a Morse function
$\tilde{f}:V_{\tilde{\Gamma}} \to \R$ by setting
\begin{equation}\label{13.1}
\tilde{f}(p,o)=f(p)
\end{equation}
and
\begin{equation}\label{13.2}
\tilde{f}(p,o)=f(p)+C,
\end{equation}
where $C$ is larger than the maximum value of $f$.  Notice that when
$$
\max(f)<c_0<C-min(f),
$$
the cross section $\tilde{\Gamma}_{c_0}$ is
just the graph $\Gamma$ itself.  Moreover, it is easy to see that the
``up'' connection on $\tilde{\Gamma}_{c_0}$ coincides with the
original connection on $\Gamma$ and the axial function on
$\tilde{\Gamma}_{c_0}$ with the original axial function. 
Now let's count
$\dim\Pi(\tilde{\Gamma}_{c_0})$ using Theorem~\ref{th:MorseLemma1}.
The critical points of $\tilde{f}$ with critical value less than $c_0$
are just the points $(p,0)$, with $p\in V_{\Gamma}$ and the index of
each of these points is simply the index of $p$.  Therefore, since
$\tilde{\Gamma}$ is a $(d+1)$-valent graph, 
the dimension of $H^r(\tilde{\Gamma}_{c})$ changes by
\begin{equation}\label{eq:A}
d_{r-k}(n+1)-d_{r-(d+1-k)}(n+1)
\end{equation}
every time one passes through a critical point of index $k$.  Thus the
total change in dimension as one goes from $c< \min\tilde{f}$ to $c_0$
is
\begin{equation}\label{eq:B}
\sum_{k=0}^d d_{r-k}(n+1)\beta_k - \sum_{k=1}^{d+1}d_{r-k}(n+1)\beta_{d+1-k}
\end{equation}
However, by Poincar\'e duality, $\beta_{d+1-k}=\beta_{k-1}$, so we can
rewrite the second sum in \eqref{eq:B} as
$$
\sum_{k=1}^{d+1}d_{r-k}(n+1)\beta_{k-1}
$$
or as
$$
\sum_{k=0}^{d}d_{r-(k-1)}(n+1)\beta_{k}.
$$
But by \eqref{eq:7K},
$$
d_{r-k}(n+1)-d_{r-(k-1)}(n+1)=d_{r-k}(n)
$$
so the combined sum \eqref{eq:B} is just
$$
\sum_{k=0}^dd_{r-k}(n)\beta_k.
$$
This completes the proof of the main theorem.
\end{proofof}

\subsection{Parallel redrawings}\label{subse:7.5}

For $r=1$, the formula \eqref{eq:7.1} tells us that $\beta_0$ is equal
to the number of connected components of $\Gamma$, a fact which is not
completely without interest.  Much more interesting, however, is the
case $r=1$.

\begin{DEF}
A {\em parallel redrawing of $\Gamma$} is a map $\pi:V\to\R^n$ such
that 
\begin{equation}\label{eq:para}
\pi(p)-\pi(q)=\lambda \alpha(p,q),\mbox{ for some } \lambda\in\R,
\end{equation}
for every edge $(p,q)\in E$. 
\end{DEF}

If $\Gamma$ is embedded graph in
$\R^n$, the identity \eqref{eq:para}
asserts that the deformation $\Gamma\mapsto \Gamma_\varepsilon$ defined by
replacing each vertex $p$ by $p+\varepsilon\pi(p)$ leaves the edges of
the deformed graph parallel to the edges of the original.

Every exact immersion is a parallel redrawing, but there are
others. In particular, the Euclidean translations are parallel
redrawings. So is dilation.

The set $\Pi(\Gamma)$ of all parallel redrawings is  
clearly a vector space, and the ``number of parallel redrawings'' is
its dimension.

Since the condition that $\pi(p)-\pi(q)$ be parallel to
$\alpha(p,q)$ is identical to the interpolation condition
$$
\pi(p)\equiv\pi(q)\mod\alpha(p,q).
$$
the dimension of $\Pi(\Gamma)$ is the 
dimension of $H^1(\Gamma,\alpha)$.  For the graphs we considered in
Theorem~\ref{th:||redrawing}, this dimension is just
\begin{equation}\label{eq:paralleldim}
n\beta_0+\beta_1.
\end{equation}

We write the count in this form even though we assumed in our proof
that $\beta_0 = 1$ since the count in that more general form is still
sometimes correct.

We can think of the first term $n\beta_0 = n$ as
counting the $n$ translations. The dilation must be a linear
combination of the other $\beta_1$ parallel redrawings.

We observed earlier that the betti numbers are preserved under
projection. So is the number of nontrivial parallel redrawings. Only
the number of translations changes, and the first term accounts for
that exactly.

\begin{figure}[h]
\centerline{
\epsfig{figure=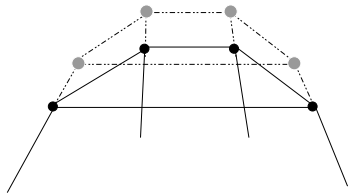,height=1in}
}
\centerline{
\parbox{4.5in}{\caption[Parallel redrawing.]{\small This
shows a parallel redrawing by moving a facet parallel to itself.}\label{fig:||}} 
}
\end{figure}

Suppose $\Gamma$ is the one-skeleton of a simple polytope $P$. Then
every facet (face of codimension $1$) determines a parallel redrawing:
just move the hyperplane containing that facet parallel to itself, as
in Figure~\ref{fig:||}.

Formula 1.1, proved in Section 5, counts the faces of $P$. When
$k=1$ it tells us that $P$ has $n\beta_0 + \beta_1$ facets, so our
construction has found all the parallel redrawings. 

\subsection{Morse inequalities}\label{subse:morseineq}

The identites \eqref{eq:7.1} imply trivially that
\begin{equation}\label{eq:ineq}
\dim H^r(\Gamma,\alpha)\leq \sum_{\ell=0}^r \dim S^{r-\ell}\beta_\ell.
\end{equation}
In the next section, we will study an analogue of $H^r(\Gamma,\alpha)$
for which these inequalities hold, but in general the equalities do
not.  We will give a new proof of these inequalities here.
This proof is much simpler
than the proof of \eqref{eq:7.1}, and in fact \eqref{eq:7.1} is true
under somewhat weaker hypotheses.
\begin{theorem}\label{th:morseineq}
Suppose that $(\Gamma,\alpha)$ is a graph with a two-independent axial
function, and that $\Gamma$ admits a $\xi$-compatibile Morse function,
$f:V\to \R$.  Then the inequalities \eqref{eq:ineq} hold.
\end{theorem}

\begin{proof}
Let $H^r(\Gamma,c)$ be the set of all interpolation schemes $<g_p>$ of
degree $r$ for which $g_p=0$ whenever $f(p)<c$.  As in
Section~\ref{subse:7.2}, let $c=f(p)$ be a critical value of index $k$
and let $c_\pm=c\pm\varepsilon$ for $\varepsilon$ small.  We claim
that there is a short exact sequence
\begin{equation}\label{eq:ses}
0\to H^r(\Gamma,c_+)\to H^r(\Gamma,c_-)\to
S^{r-k}\alpha_{e_1}\cdots\alpha_{e_k}
\end{equation}
where $e_i$, $i=1,\dots,k$, are the $k$ edges in $\star (p)$ with
$\alpha(e_i)\cdot \xi<0$.   It is clear that $H^r(\Gamma,c_+)$ is the
kernel of the map
$$
<g>\in H^r(\Gamma,c_-)\mapsto g_p \in S^r.
$$
But if $g$ is in $H^r(\Gamma,c_-)$ then for every down-pointing edges
$e_i$ with terminal vertex $q_i$, $f(q_i)<c_-$, so $g_{q_i}=0$.
Hence, the interpolation conditions
$$
g_p\equiv g_{q_i}\mod\alpha(e_i)
$$
imply that $g_p$ is divisible by $\alpha(e_i)$.  Moreover, since the
$\alpha(e_i)$'s are two-independent, $g_p$ is divisible by
$\alpha(e_1)\cdot\cdots\cdot\alpha(e_k)$, completing the proof that
the sequences \eqref{eq:ses} is exact.

This short exact sequence implies
$$
\dim H^r(\Gamma,c_-)-\dim H^r(\Gamma,c_+)\leq\dim S^{r-k} \ .
$$
If $c_0<\min(f)<\max(f)<c_1$, summing these
inequalities yields the desired 
$$
\dim H^r(\Gamma,\alpha)=\dim H^r(\Gamma,c_0)-\dim H^r(\Gamma,c_1)\leq
\sum\dim(S^r-k)\beta_k \ .
$$
\end{proof}

\section{Interpolation schemes involving polynomials in
``anti-commuting'' variables} \label{sec:8}

The interpolation schemes that we considered in Sections~\ref{se:poly} and
\ref{se:para} involved polynomials, $f \in S^r$ of the  form
\begin{displaymath}
f= \sum_{i_1 + \cdots + i_n =r} a_{i_1 \ldots i_n}
     x_1^{i_1} \ldots x_n^{i_n}
\end{displaymath}
in ``commuting'' variables:  i.e.,~$x_ix_j=x_jx_i$.  One can also,
however, consider interpolation schemes involving polynomials in
``anti-commuting'' variables 
\begin{displaymath}
  f= \sum_{i_1 + \cdots + i_n=r} a_{i_1} \ldots a_{i_n}
     x_{i_1} \ldots x_{i_n}
\end{displaymath}
 where $x_ix_j =-x_jx_i$; in other words interpolation schemes,
$<f_p>$, for $p \in V$ in which $f_p$ sits in the $r$\st{th} exterior
power, $\Lambda^r (\R^n)$, of the vector space, $\R^n$.  The
interpolation conditions 
 \begin{equation}\label{eq:fermionic}
   f_p \equiv f_q \mod \alpha_e
 \end{equation}
still make sense in this anti-commuting context if we interpret this
equation as saying that $f_p-f_q \in \alpha_e \wedge \Lambda^{r-1}$. 

Let us denote by $\tilde{H}^r (\Gamma ,\alpha)$ the space of all
$r$\st{th} degree solutions of the equations (\ref{eq:fermionic}).  The sum 
\begin{displaymath}
  \tilde{H} (\Gamma , \alpha) = \oplus^n_{r=0}
    \tilde{H}^r (\Gamma ,\alpha)
\end{displaymath}
is, like its ``bosonic'' counterpart, $H(\Gamma ,\alpha)$, a graded
ring; and, in particular, a graded module over the exterior algebra
$\Lambda (\R^n)$.  We claim that the following analogue of
Theorem~\ref{th:morseineq} is true. 

\begin{theorem}
  \label{th:8}
  Suppose that $\alpha$ is $r$-independent and that there exists an
injective $\xi$-compatible Morse function, $f:V \to \R$.  Then
  \begin{equation}
    \label{eq:8.1}
    \dim \tilde{H}^r (\Gamma ,\alpha) \leq \sum^r_{k=0}
    \binom{n-k}{n-r} b_k \, .
  \end{equation}

\end{theorem}

\begin{proof}
The ``fermionic'' analogue of (\ref{eq:ses}) asserts that there is an
exact sequence
\begin{displaymath}
  0 \to \tilde{H}^r (\Gamma , C_+) \to \tilde{H}^r (\Gamma ,C_-)
     \to \Lambda^{r-k} \alpha_{e_1} \wedge \ldots \wedge \alpha_{e_k} \, .
\end{displaymath}
Hence  
\begin{displaymath}
 \dim \tilde{H}^r (\Gamma , C_-) - \dim \tilde{H}^r (\Gamma ,C_+)
    \leq \binom{n-k}{n-r}
\end{displaymath}
and by adding up these identities as in Section~\ref{subse:morseineq}
one gets (\ref{eq:8.1}).
\end{proof}

Are these Morse inequalities ever equalities?  We will show that they
are if $\Gamma$ is the one-skeleton of a simple convex $n$-dimensional
polytope, $\Delta$; and, in fact, we will show in this case that these
identities are identical with the McMullen identities (\ref{eq:1}).
The idea of our proof will be to show that for every $n-r$ dimensional
face, $F$, one can associate an $r$\st{th} order solution of the
interpolation identities (\ref{eq:6.1}) by mimicking the constructing in
Section~\ref{se:6.4}, and by showing that the solutions constructed this
way form a basis of $\tilde{H}^r (\Gamma ,\alpha)$.  Fix a set of
vectors, $v_1 , \ldots ,v_r \in \R^n$ such that the $v_r$ 's are
normal to the face $F$ and are linearly independent; and for every
vertex, $p$, of $\Gamma$ define
\begin{equation}
  \label{eq:8.2}
  f^F_p =0
\end{equation}
if $p$ is not a vertex of $F$ and
\begin{equation}
  \label{eq:8.3}
  f^F_p = c_p \alpha_{e_1} \wedge \ldots \wedge \alpha_{e_r}
\end{equation}
if $p$ is a vertex of $F$, the $e_i$'s being, as in
Section~\ref{se:6.4}, the edges of $\Gamma$ normal to $F$ at $p$ and
$c_p$ being defined by the normalization condition
\begin{equation}
  \label{eq:8.4}
  c_p \det (\alpha_{e_i} (v_j)) =1 \, .
\end{equation}
It is easy to check that (\ref{eq:8.2})--(\ref{eq:8.4}) is a solution
of the interpolation equations (\ref{eq:8.1}).  We will prove the
following theorem.

\begin{theorem}
Let $F_i , i=1,\ldots ,N$ be the $(n-r)$-dimensional faces of
$\Lambda$.  Then the interpolation schemes, $f^{F_i}, i=1,\ldots ,N$
are a basis of $\tilde{H}^r (\Gamma ,\alpha)$. 
\end{theorem}

\begin{proof}
By McMullen's identity the right-hand side of (\ref{eq:8.1}) is equal
to $N$ so it suffices to prove that the $f^{F_i}$'s are linearly
independent; and hence it suffices to prove that at each vertex, $p$,
the vectors
\begin{displaymath}
  f^{F_i}(p)\, , \quad p \hbox{  a vertex of  } F_i
\end{displaymath}
are linearly independent.  However 
it is clear that these vectors are in fact a basis of $\Lambda^r (\R^n)$.
\end{proof}

\section{Examples}\label{se:EGs}

In this final section we review the examples we have been following
through the text and introduce some new ones that suggest new
directions to explore. When proofs are short we include them.
Some will be found in \cite{Hthesis}. Others we leave as exercises.

\subsection{The complete graph $K_n$.}

Our standard view $K_n$ embeds with vertices the standard
basis vectors in $\R^{n}$. That embedding is a regular simplex in the
$n-1$-dimensional subspace $\Sigma x_i = 1$. The exact axial function
is determined by assigning to each vertex the difference between its
end points.  The following figure shows a part of the connection
determined by that axial function for $K_4$: it moves edges across the
triangular faces.

\begin{figure}[h]
\centerline{
\psfig{figure=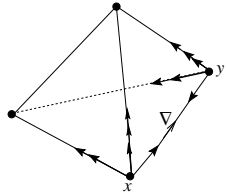,height=1in}
}
\centerline{
\parbox{4.5in}{\caption[Parallel Transport on $K_4$.]{{\small This
shows the connection we defined above on the graph $K_4$.}}}
}
\end{figure}
When we think of $K_4$ just as an abstract $4$-regular graph we 
find that it has 10 different connections (up to graph
automorphism). But in each of these connections other than the
standard one there is at least one geodesic of length at least $4$, so
none of those connections has a $3$-independent immersion. So we will
study only the standard view.

\begin{prop}
\label{geo}
The geodesics of $K_n$ are the triangles. The connected totally
geodesic subgraphs are the complete subgraphs.
\end{prop}

\begin{proof}
It's clear that the geodesics are the triangles.
Let $\Gamma_0$ be a connected totally geodesic subgraph
and $p$ and $q$ two vertices of $\Gamma_0$. Then transporting edge
$e = (p,q)$ along a path in $\Gamma_0$ from $p$ to $q$ we eventually reach
a triangle containing $q$. At that point the image of $e$ transports
to an edge of $\Gamma_0$ so $e$ must have been part of $\Gamma_0$ to
begin with.
\end{proof}

\noindent It's easy to compute the holonomy of $K_n$.

\begin{prop}
$Hol(K_n)\cong S_{n-1}$. 
\end{prop}

\begin{proof}
If you follow the connection along triangle $(p,q,r)$ from $p$ back to
itself you interchange $(p,q)$ and $(p,r)$. Thus the holonomy group
acting on $\star(p)$ contains all the transpositions.
\end{proof}

\begin{prop}
The Betti numbers of $K_n$ are invariant of choice of direction $\xi$
and are $(1,1,\dots,1)$.
\end{prop}

\begin{proof}
The geodesics are triangles, hence convex. hence inflection free, so
the Betti numbers are well defined. 
Let $\xi=(1, 2, \dots, n)$.  Then the number of down edges at the vertex
corresponding to the $i^{th}$ coordinate vector is the
number of $j$'s less than $i$.
\end{proof}

\subsection{The Johnson graph $J(n,k)$}

Recall that the Johnson graph $J(n,k)$ is the graph with vertices
corresponding to $k$-element subsets of $\{1,2,\dots ,n\}$; two
vertices $S,T\in V$ are adjacent if \/ $\# (S\cap T)=k-1$. Then we can
think of an oriented edge as an ordered pair $(i,j)$: to get
from $S$ to $T$ we remove $i$ and add $j$. We naturally embed
$J(n,k)$  in  $\R^{n}$ by representing each vertex as a vector with
$k$ $1$'s and $n-k$ $0$'s. 
That embedding is a $k \times (n-k)$-regular polytope in the
$n-1$-dimensional subspace $\Sigma x_i = k$. The exact axial function
is determined by assigning to each vertex the difference between its
end points. 

The easiest way to describe the natural connection is to describe its
geodesics. They are the triangles $Q \cup \{a\}$, $Q \cup \{b\}$,
$Q \cup \{c\}$ for $k-1$  element sets $Q$ and distinct $a$, $b$,
$c$ and the planar squares $Q \cup \{a,b\}$, 
$Q \cup \{b,c\}$, 
$Q \cup \{c,d\}$, 
$Q \cup \{d,a\}$, 
for $k-2$  element sets $Q$ and distinct $a$, $b$,
$c$, $d$.

The triangles are actual faces of the polytope. The squares are more
like equators, as in the picture of the octahedron $J(4,2)$ below.

\begin{figure}[h]
\centerline{
\epsfig{figure=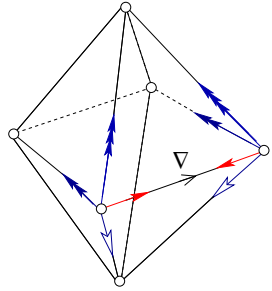,height=1.5in}
}
\centerline{
\parbox{4.5in}{\caption[Parallel Transport on Gr(4,2).]{{\small This
shows the connection we defined on the Johnson graph $J(4,2)$.}}} 
}
\end{figure}
Knowing this connection it is not hard to find all 
the totally geodesic subgraphs of $J(n,k)$.

\begin{prop}
If $\Gamma_0$ is a totally geodesic subgraph of $\Gamma=J(n,k)$, then
$$
\Gamma_0\cong J(A_1,\ell_1)\times\cdots\times J(A_r,\ell_r),
$$
where the $A_i$ are subsets of $\{ 1,\dots ,n\}$ of size
$a_i\geq\ell_i$, and $\{ 1,\dots ,n\}$ is the disjoint union of the
$A_i$.
\end{prop}

We can also compute the holonomy of $J(n,k)$:

\begin{prop}
$Hol(J(n,k))\cong S_k\times S_{n-k}$.
\end{prop}

The proof is similar to that for the complete graph.

Since the geodesics are convex the Betti numbers are well defined. For
the octahedron they are $(1,1,2,1,1)$. Since $n\beta_0 + \beta_1 = 3
\times 1 + 1 = 4$ the octahedron has only trivial parallel redrawings.

We leave the computation of the Betti numbers of the general Johnson
graph as an exercise for the reader.

\subsection{The permutahedron $S_n$.}

The permutahedron is the Cayley graph of the symmetric
group $S_n$, generated by reflections. Its vertices 
correspond to permutations of $\{1, \ldots , n \}$. Two
vertices (permutations) are adjacent if they differ by left
multiplication by a transposition $(i,j)$: the vertices joined
to $\sigma \in S_n$ are the permutations $\tau \cdot \sigma$ for all
transpositions $\tau$. This construction defines a natural labeling of
the edges in $\star(\sigma)$ by the
transpositions. In turn that defines a natural connection: map an edge
in the star of a vertex to the edge labelled by the same transposition
in the star of an adjacent vertex.

To embed $S_n$ we think of its vertices as permutations of
the entries of the vector $(1, \ldots , n )$.
The convex hull of
that embedding is a simple 
polytope in the
$n-1$-dimensional subspace $\Sigma x_i = n(n+1)/2$, but the convex
hull is only a small part of the story. The graph is
$\binom{n}{2}$-regular. Most of the embedded edges are internal to the
polytope. When so embedded the natural axial function is exact, with
inflection free geodesics. 

\begin{figure}[h]
\centerline{
\epsfig{figure=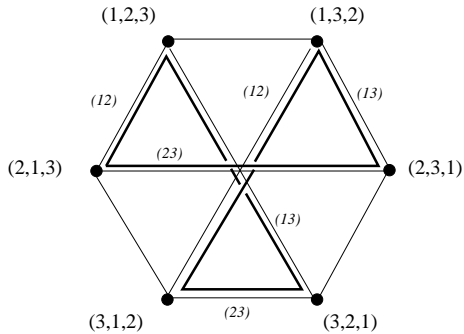,height=1.75in}
}
\centerline{
\parbox{4.5in}{\caption[A geodesic on $S_3$.]{\small 
An immersion of $S_3$, with one of the three geodesics in bold. The
edges are labelled by transpositions.}\label{fig:XXXX}}
}
\end{figure}

$S_3$ is a regular hexagon together with its diagonals. Figure~\ref{fig:XXXX}
shows $S_3$ and one of its three geodesics.  Since each of the
geodesics is inflection free (although not convex!) the Betti numbers
are well defined. They are $(1,2,2,1)$.  In the next section we will
discuss the parallel redrawings of $S_3$.

The permutahedron $S_4$ is the truncated octahedron in $R^3$ together
with the necessary internal edges. It has hexagonal and square faces
(internal and external) corresponding to natural subgroups $S_3 \times
S_1$ and $S_2 \times S_2$ of $S_4$. 
These are totally geodesic subgraphs. The reader can decide whether
the analagous construction produces all totally geodesic subgraphs.
$S_4$ has Betti numbers
$(1,3,5,6,5,3,1)$. In general

\begin{prop}
The generating polynomial $B_n(z)$ for the Betti numbers of $S_n$ is
$$
	B_n(z) = (1 + z + \cdots + z^{n-1}) B_{n-1}(z)
$$

$$
	= \prod_{k=0}^{n-1}(1 + z + \cdots + z^{k})
$$
\end{prop}
\begin{proof}
If we compute the Betti numbers for $S_n$ using the generic direction 
$\xi = (1,2, \ldots , n)$ then the number of down edges at $\sigma$ is
the number of inversions in $\sigma$. The generating function in the
proposition is the one that counts permutations according to the
number of inversions. \cite{sloane}
\end{proof}

The holonomy of $S_n$ is trivial:

\begin{prop}
$Hol(S_n)\cong 0$.
\end{prop}
\begin{proof}
The connection just matches edges labelled by the same
transposition, so following a chain from a vertex back to itself
permutes nothing in the star of that vertex.
\end{proof}

The permutahedra are examples from the class
of Cayley graphs which have
a GKM graph structure which is compatible with their structure as a
Cayley graph. $S_n$ corresponds to the full flag manifold of all  
subspaces of $CP^n$.
Cayley graphs are discussed more thoroghly in
\cite{Hthesis}.

\subsection{The complete bipartite graph $K_{n,n}$.}

Let $\mathcal D = \D$ be the group of symmetries of the regular
$n$-gon: the dihedral
group with $2n$ elements. Then $\D$ is a reflection group of type
$I_2(n)$, following the notational conventions of Humphreys
\cite{Hu}. It is generated by two reflections, and contains $n$
reflections and $n$ rotations.  If we let $\Delta$ be the set of
reflections in $\D$, then the Cayley graph $\Gamma=(\D,\Delta)$
has vertices corresponding to elements of $\D$. $\sigma \in \D$
is connected to $ \tau\sigma$ for every $\tau \in \Delta$.
Just half the vertices of $\Gamma$ correspond to symmetries that preserve
the orientation of the $n$-gon, and 
$\sigma$ preserves 
orientation if and only if $\tau\sigma$ reverses it. Thus
the graph is
bipartite.  The only $n$-regular bipartite graph on $2n$ vertices is
$K_{n,n}$.  

$D_n$ has a natural holonomy free connection defined just as for the
permutahedron, using the reflection generating one vertex from
another as the label for the corresponding edge. The natural embedding
of $D_n$ as the vertices of a regular 
$2n$-gon produces an exact axial function with
inflection free geodesics for that connection.

$\mathcal{D}_3$ is $K_{3,3}$ and also the permutahedron $S_3$ discussed above.
The figure below shows two more examples.

This class of graphs is particularly interesting because $\mathcal{D}_n =
K_{n,n}$ is the graph associated with a manifold 
in the sense described in Section~\ref{subse:GKM} only when
$n = 1,2,3,4,6$, so they provide examples where combinatorics may go
further than differential geometry.  

\begin{figure}[h]
\centerline{
\epsfig{figure=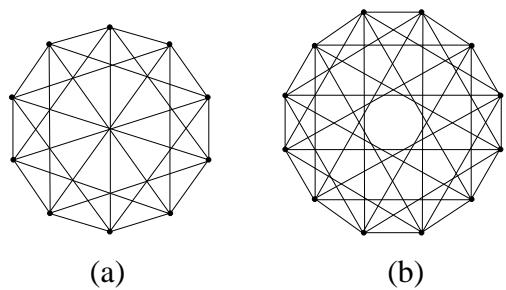,height=1.5in}
}
\centerline{
\parbox{4.5in}{\caption[Parallel Transport on Gr(4,2).]{\small This
shows the Cayley graphs for (a) $\mathcal{D}_5$ and (b) $\mathcal{D}_6$.}\label{fig?}} 
}
\end{figure}

We will leave as an exercise the following Betti number count.

\begin{prop}
The Betti numbers of $K_{n,n}$ are invariant of choice of direction $\xi$
and are $(1,2,\dots,2,1)$.
\end{prop}

\begin{figure}[h]
\centerline{
\epsfig{figure=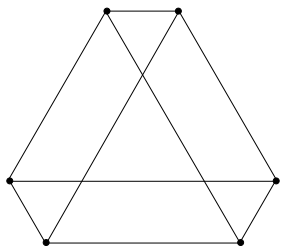,height=1.25in}
}
\centerline{
\parbox{4.5in}{\caption{\small A parallel redrawing of $K_{3,3}$.}\label{fig:YYY}}
}
\end{figure}

Note that $\D$ is far from $3$-independent.
Nevertheless $n\beta_0 + \beta_1 = 2 \times 2 + 2 = 4$ counts the
number of parallel redrawings. There are the three trivial ones (two
translations and the dilation) and one significant one:
rotate the sense preserving symmetries clockwise and others
counterclockwise 
around the circle on which they lie. Figure~\ref{fig:YYY} shows the resulting
deformation of the hexagon.

Explaining the existence of this deformation is worth a
digression. The first remark is that the hypothesis of
$3$-independence in Theorem~\ref{th:||redrawing} can be dropped if we
assume instead that the axial function is exact, as it is in these
cases. For a proof in the context of GKM manifolds, see \cite{GZ2}. Thus we
do expect to see a nontrivial parallel deformation here.

The theory of rigidity predicts the same deformation. An 
{\em infinitesimal motion} of an embedded graph is an assignment of
a velocity vector to each vertex in such a way that the length of each
edge is (infinitesimally) unchanged. The space of infinitesimal
motions includes the Euclidean motions and perhaps others.
For precise definitions see \cite{AR}. In the plane there
is a one to one correspondence between infinitesimal motions and
parallel redrawings: rotating all the vectors of an infinitesimal
motion through a quarter of a turn converts that motion into a
parallel redrawing. Translations remain translations. Rotation
becomes dilation. Nontrivial motions become nontrivial parallel 
redrawings.

A nineteenth century
theorem (reproved and generalized in \cite{BR}) says that a plane
embedding of  $K(m,n)$
(for $m, n \ge 3$)
is rigid (no infinitesimal motions) except when it lies
on a conic. In this case that's just what happens. The regular $2n$-gon
lies on a circle. 
The single nontrivial infinitesimal motion moves the odd
permutations radially outward while the even ones move inward at the
same velocity. Rotating that motion a quarter of a turn produces the
parallel redrawing: half the vertices move clockwise, half
counterclockwise. Figure~\ref{fig:YYY} shows the result for $S_3$.

Finally, even without the picture we could have deduced the existence
of the conic on which $S_3$ lies from the the exactness of the axial
function together with the dual of Pascal's theorem in projective
geometry.

\subsection{Several more examples}

We conclude our bestiary with several final suggestive examples.

Whenever the plane containing two adjacent edges of a polytope
intersects that polytope in a cycle of edges the one skeleton of the
polytope is an embedded graph with a natural axial function. The
polytope need not be simple. The cuboctahedron, shown below, provides one
example. It is $4$-regular with $6$ square and $8$ triagular faces.
Three hexagonal plane sections define three more geodesics. 
Its Betti numbers are $(1,2,6,2,1)$ so it has one nontrivial parallel
redrawing, which dilates four of the triangular faces, converting the
square faces to rectangles. 

\begin{figure}[h]
\centerline{
\psfig{figure=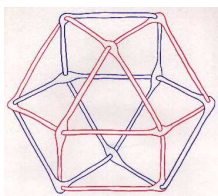,height=1in}
}
\centerline{
\parbox{4.5in}{\caption[Not stellated icos]{\small
The cuboctahedron.}} 
}
\end{figure}

\begin{figure}[h]
\centerline{
\epsfig{figure=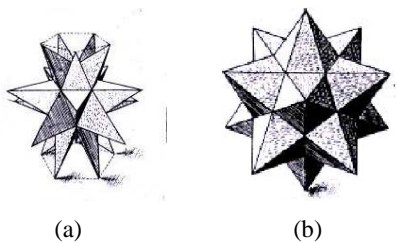,height=1.25in}
}
\centerline{
\parbox{4.5in}{\caption[Not stellated icos]{\small
This shows (a) the great stellated dodecahedron (a stellated icosahedron) and (b) small
stellated dodecahedron (a stellated dodecahedron).}\label{figK}} 
}
\end{figure}

Figure~\ref{figK}(a) shows 
the {\em great stellated dodecahedron}, from Kepler's 
1619 {\em Harmonice Mundi}. 
It is in fact a stellated icosahedron. It's 
a simple polytope with pentagrams for faces. These are the
geodesics. Since these are inflection free the
Betti numbers are well defined. They are $(5, 5, 5, 5)$. In this case 
$n\beta_0 + \beta_1 = 3\cdot 5+5=20$ does properly count the
number of parallel redrawings, since there is one for each of the $20$
faces. However, Theorem~\ref{th:||redrawing}
does not apply because the zeroth Betti number for each geodesic is
$2$, not $1$.

Figure~\ref{figK}(b) shows
Kepler's {\em small stellated dodecahedron}, which {\em is} a stellated
dodecahedron,
Its geodesics are $12$ pentagrams and $20$ equilateral triangles. The
invariant Betti numbers are $(3,1,2,2,1,3)$. In this case $n\beta_0 +
\beta_1 = 3 \cdot 3 + 1 = 10$ counts neither the number of faces in
Kepler's sense (as it would for a simple polytope) nor the number of
parallel redrawings.

In the plane, however, we can often get a correct count even when
hypotheses fail, Any $n$-gon has $n$ parallel redrawings, one for each
edge, and $n = 2\beta_0 + \beta_1$ as long as Poincare duality holds,
even when the Betti numbers are not invariant. Figure~\ref{fig:ZZZ} below shows
the dart, whose Betti numbers are $(1,2,1)$ or $(2,0,2)$ depending on
the choice of $\xi$. 

\begin{figure}[h]
\centerline{
\epsfig{figure=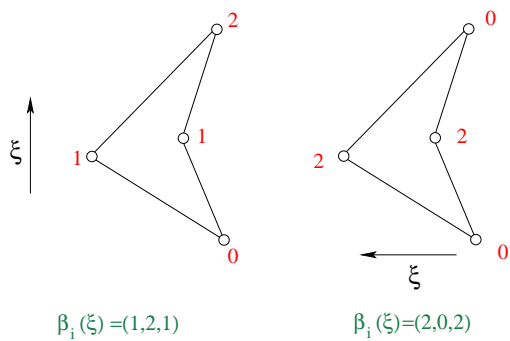,height=1.75in}
}
\centerline{
\parbox{4.5in}{\caption[Not stellated icos]{\small
This shows the Betti numbers for the dart, with its Betti numbers
$\beta_i(\xi)$ for two choices of $\xi$.  The number at each vertex is
the index, and the Betti numbers are given below each figure.}\label{fig:ZZZ}} 
}
\end{figure}

\begin{figure}[h]
\centerline{
\epsfig{figure=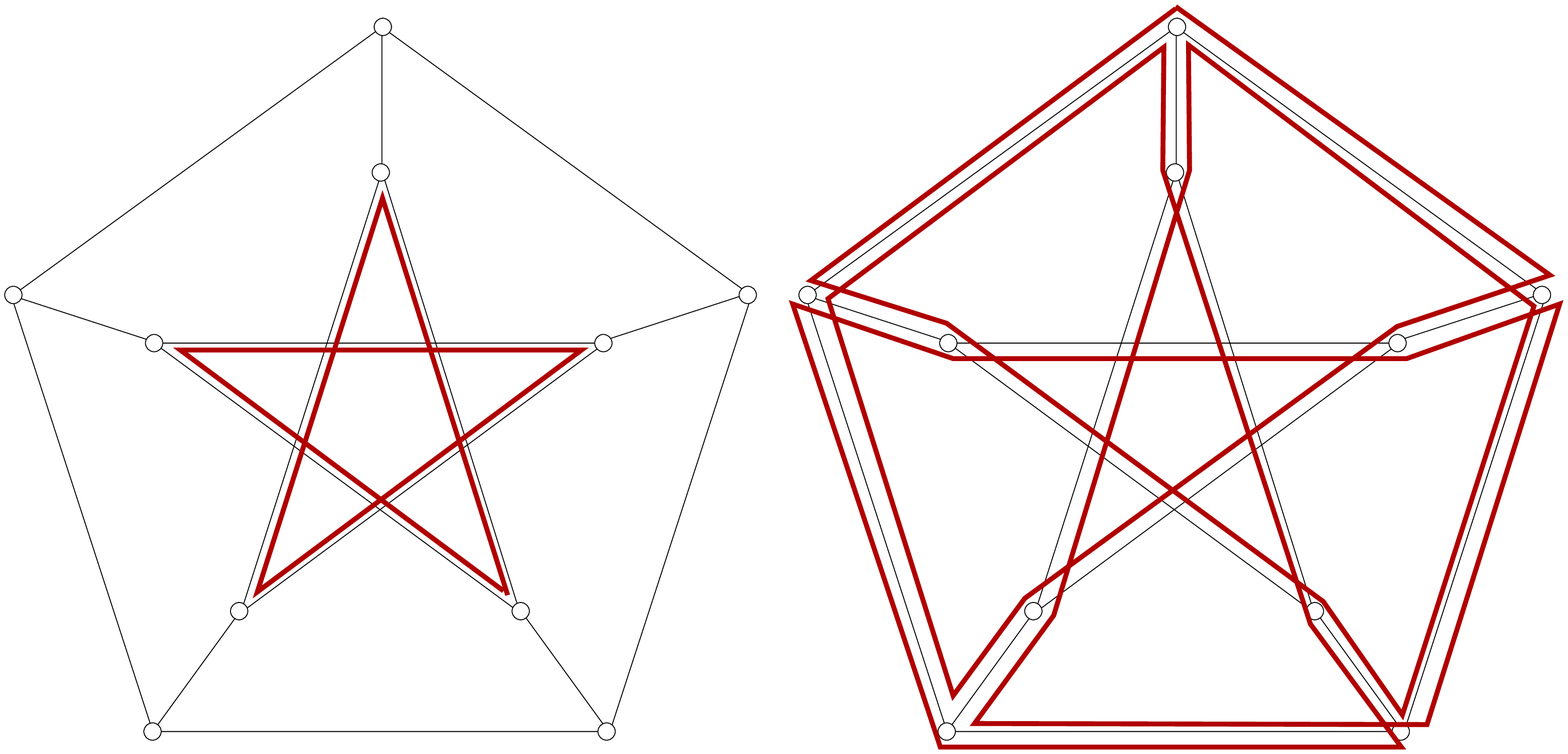,height=1in}
}
\centerline{
\parbox{4.5in}{\caption[Petersen graph]{\small The
two geodesics of the Petersen graph.}\label{figP}} 
}
\end{figure}

Finally
Figure~\ref{figP} shows the 
Petersen graph in the plane with two inflection free geodesics that
define a connection. 
The well defined Betti numbers are
$(1,4,4,1)$. Since the axial 
function is exact we can count parallel redrawings even though it is
not $3$-independent. There are $2 \cdot 1 + 4 = 6$. Five correspond
to edges of the enclosing pentagon which, when moved independently,
force parallel redrawing of the inner pentagram. The sixth is a
dilation of the inner pentagram while the outer pentagon remains fixed.
It corresponds to the infinitesimal motion that rotates the inner
pentagram relative 
to the outer pentagon - an infinitesimal motion possible only because
exactness means the radial edges will meet in a point when extended.

\end{document}